\documentclass[12pt, twoside, a4paper]{amsart}

\usepackage{latexsym, amsfonts, amsmath, amscd, amssymb, amsthm, epsfig, graphics}

\newtheorem{thm}{Theorem}[section]
\newtheorem{prop}[thm]{Proposition}
\newtheorem{lem}[thm]{Lemma}
\newtheorem{cor}[thm]{Corollary}

\theoremstyle{definition}
\newtheorem{defn}[thm]{Definition}
\newtheorem{ex}[thm]{Example}
\newtheorem{rem}[thm]{Remark}

\newcommand{\thmref}[1]{Theorem~\ref{#1}}
\newcommand{\lemref}[1]{Lemma~\ref{#1}}
\newcommand{\propref}[1]{Proposition~\ref{#1}}
\newcommand{\defref}[1]{Definition~\ref{#1}}
\newcommand{\exref}[1]{Example~\ref{#1}}
\newcommand{\remref}[1]{Remark~\ref{#1}}
\newcommand{\corref}[1]{Corollary~\ref{#1}}
\newcommand{\figref}[1]{Figure~\ref{#1}}
\newcommand{\secref}[1]{Section~\ref{#1}}

\newcommand{\lra}{\longrightarrow}
\newcommand{\lmt}{\longmapsto}

\newcommand{\CC}{\mathcal{ C}}
\newcommand{\DD}{\mathcal{ D}}
\newcommand{\EE}{\mathcal{ E}}
\newcommand{\FF}{\mathcal{ F}}
\newcommand{\LL}{\mathcal{ L}}

\newcommand{\WW}{\mathcal{ W}}

\newcommand{\al}{\alpha}

\newcommand{\R}{\mathbb{ R}}

\newcommand{\N}{\mathbb{ N}}
\newcommand{\Z}{\mathbb{ Z}}
\newcommand{\PP}{\mathbb{ P}}

\newcommand{\eing}[1]{\big|_{#1}}
\newcommand{\rve}[1]{\frac{\partial}{\partial{#1}}}
\newcommand{\ww}{\wedge}
\newcommand{\x}{\times}
\newcommand{\eps}{\varepsilon}
\newcommand{\SO}{\mathrm{SO}}

\DeclareMathOperator{\pr}{pr}
\DeclareMathOperator{\id}{id}
\DeclareMathOperator{\tb}{tb}

\DeclareMathOperator{\tw}{tw}
\DeclareMathOperator{\fr}{fr}
\DeclareMathOperator{\rot}{rot}

\newcommand{\mlabel}[1]{%\marginpar{#1} 
                        \label{#1}}

\title{Existence of Engel structures}
\author{Thomas Vogel}
\email{tvogel@mathematik.uni-muenchen.de}

\begin{document}

\begin{abstract}
We develop a construction of Engel structures on $4$--manifolds based on decompositions of manifolds into round handles. This allows us to show that all parallelizable $4$--manifolds admit an Engel structure. We also show that, given two Engel manifolds $M_1,M_2$ satisfying a certain condition on the characteristic foliations, there is an Engel structure on $M_1\# M_2\# (S^2\x S^2)$ which is closely related to the original Engel structures.
\end{abstract}

\maketitle

\tableofcontents

\section{Introduction}
Distributions are subbundles of the tangent bundle of a manifold. It is natural not to consider general distributions but to make geometric assumptions, for example integrability. In this case the distribution defines a foliation on the manifold. Another possibility is to assume that a distribution is nowhere integrable. Important examples of this type are contact structures on odd--dimensional manifolds. On $3$--dimensional manifolds, properties of contact structures reflect topological features of the underlying manifold in a surprising way. 

Engel structures form another class of non--integrable distributions which is closely related to contact structures. By definition, an Engel structure is a smooth distribution $\DD$ of rank $2$ on a manifold $M$ of dimension $4$ which satisfies the non--integrability conditions
\begin{align*}
\mathrm{rank}[\DD,\DD]&=3 \\
\mathrm{rank}[\DD,[\DD,\DD]]&=4~, 
\end{align*}
where $[\DD,\DD]$ consists of those tangent vectors which can be obtained by taking commutators of local sections of $\DD$. Examples of Engel structures can be obtained from contact structures on $3$--manifolds. In this article we use results from contact topology to investigate Engel structures. This shows that Engel structures and contact structures are closely related.

Similarly to contact structures and symplectic structures, all Engel structures are locally isomorphic. Every point of an Engel manifold has a neighborhood with local coordinates $x,y,z,w$ such that the Engel structure is the intersection of the kernels of the $1$--forms
\begin{align} \label{e:engel normal}
\begin{split}
\alpha & = dz-x\,dy \\
\beta  & = dx-w\,dy~.
\end{split}
\end{align}
This normal form was obtained first by F.~Engel in \cite{engel}. Together with the fact that a $C^2$--small perturbation of an Engel structure is again an Engel structure, this implies that Engel structures are stable in the sense of singularity theory. In \cite{mo1} R.~Montgomery classified all distribution types with this stability property. They belong to one of the following classes:
\begin{itemize}
\item foliations of rank one on manifolds of arbitrary dimension
\item contact structures on manifolds of odd dimension
\item even contact structures on manifolds of even dimension
\item Engel structures on manifolds of dimension $4$.
\end{itemize}
Thus Engel structures are special among general distributions and even among the stable distribution types they seem to be exceptional since they are a peculiarity of dimension $4$. 

On the other hand Engel structures appear quite naturally. For example, a generic plane field on a $4$--manifold satisfies the Engel conditions almost everywhere and one can construct an Engel structure from a contact structure on a $3$--manifold. Certain non--holonomic constraints studied in classical mechanics also lead to Engel structures.

An Engel structure $\DD$ induces a distribution of hyperplanes $\EE=[\DD,\DD]$ which is an even contact structure, i.e. $[\EE,\EE]$ is the whole tangent bundle. Moreover, to the even contact structure $\EE$ one can associate its characteristic foliation $\WW$ and it turns out that $\WW$ is tangent to $\DD$. This is explained in \secref{s:engeldef}. Thus an Engel structure $\DD$ on $M$ induces a flag of distributions 
\begin{equation} \label{e:flag1}
\WW\subset\DD\subset\EE=[\DD,\DD]\subset TM 
\end{equation}  
such that each distribution has corank one in the next one. Hypersurfaces transverse to $\WW$ carry a contact structure together with a Legendrian line field. We will refer to this line field as the intersection line field. 

The presence of the distributions in \eqref{e:flag1} leads to strong restrictions for the topology of Engel manifolds. The following result can be found in \cite{montgomery2} where it is attributed to V.~Gersh\-ko\-vich.
\begin{prop} \mlabel{p:trivial intro}
If an orientable $4$--manifold admits an orientable Engel structure, then it has trivial tangent bundle.
\end{prop}
We develop a construction of Engel manifolds which allows us to prove the converse of \propref{p:trivial intro}. 
\begin{thm} \mlabel{t:exist intro}
Every parallelizable $4$--manifold admits an orientable Engel structure. 
\end{thm}
This solves a problem from \cite{mishachev eliashberg} and \cite{montgomery2}. For open manifolds, \thmref{t:exist intro} can be proved using the $h$--principle for open, $\textrm{Diff}$--invariant relations, cf.~\cite{mishachev eliashberg}. 

In the literature one can find only two other constructions of Engel structures on closed manifolds. The first one is called prolongation and was introduced by E.~Cartan. With this method one gets Engel structures on certain $S^1$--bundles over $3$--dimensional contact manifolds. The second construction is due to H.--J.~Geiges, cf.~\cite{geiges}. It yields Engel structures on parallelizable mapping tori. These two constructions cover only a small portion of all parallelizable $4$--manifolds. 

Our construction is based on decompositions of manifolds into round handles. A round handle of dimension $n$ and index $l\in\{0,\ldots,n-1\}$ is 
$$
R_l=D^l\x D^{n-l-1}\x S^1~.
$$
It is attached to a manifold with boundary using an embedding of $\partial_-R_l = S^{l-1}\x D^{n-l-1}\x S^1$ into the boundary of $M$. y a theorem of D.~Asimov in \cite{asimov1} every $4$--manifold with vanishing Euler characteristic can be obtained from $R_0$ by attaching round handles of higher index successively.

We fix a set of model Engel structures on round handles such that the characteristic foliation is transverse to $\partial_-R_l$ and $\partial_+R_l=\partial R_l\setminus\partial_-R_l$. These model Engel structures are constructed using a perturbed version of the prolongation construction. The contact structure on each $\partial_-R_l$ depends essentially only on $l$ while the intersection line field and the orientation of the contact structure are different for different model Engel structures on $R_l$.

When we attach a round handle $R_l$ with a model Engel structure to an Engel manifold $M'$ with transverse boundary we have to ensure that the attaching map preserves the oriented contact structures and the intersection line field if we want to extend the Engel structure from $M'$ to $M'\cup R_l$ by the model Engel structure. In order to satisfy these conditions we isotope the attaching map and we choose a model Engel structure suitably. For this we use several constructions of contact topology. It turns out that it is convenient to ensure that the contact structures on boundaries transverse to the characteristic foliation are overtwisted throughout the construction.

Since $M'\cup R_l$ is again an Engel manifold with transverse boundary we can iterate this construction. For the construction of an Engel structure on a manifold $M$ with trivial tangent bundle we fix a round handle decomposition of $M$ and we use the model Engel structures to extend the Engel structure when the round handles are attached successively. The condition that $M$ has trivial tangent bundle will be used to show that there is a model Engel structure with the desired properties in our collection of model Engel structures on $R_l$ for $l=1,2,3$.      

Another result which can be obtained using round handles is the following theorem. 
\begin{thm} \mlabel{t:M+N intro}
Let $M_1, M_2$ be manifolds with Engel structures $\DD_1,\DD_2$ such that the  characteristic foliations admit closed transversals $N_1, N_2$. Then $M_1\#M_2\#(S^2\x S^2)$ admits an Engel structure which coincides with $\DD_1$ and $\DD_2$ away from a neighborhood of $N_1, N_2$ where all connected sums are performed. The characteristic foliation of the new Engel structure again admits a closed transversal.
\end{thm}
This theorem can be used to construct Engel manifolds -- like for example $(n+1)(S^3\x S^1)\# n(S^2\x S^2)$ for $n\ge 1$ -- which are not covered by the construction of Geiges or by prolongation. 

This article is organized as follows. \secref{s:defengel} contains several definitions and properties of Engel structures and related distributions. In \secref{ss:vert mod} we introduce vertical modifications of the boundary of Engel manifolds. This construction will be used frequently in later sections. Our construction relies on several facts from $3$--dimensional contact topology. We summarize them in \secref{s:contact top} for the convenience of the reader. The theorems discussed here will be used later in order to bring the attaching maps of round handles in a suitable form. 

In \secref{s:round model} we discuss round handle decompositions of manifolds and model Engel structures on round handles. The various isotopies of attaching maps for round handles and the choice of the right model Engel structure are explained in \secref{s:isotopy+choice}. We give the proofs of \thmref{t:exist intro} and \thmref{t:M+N intro} in \secref{s:proofs}. 

{\em Acknowledgements: }
This article contains the main results of the authors thesis. It is a pleasure for me to thank my advisor Dieter Kotschick for his continuous support and help. I would like to thank Kai Cieliebak, Yakov Eliashberg and Paolo Ghiggini for helpful discussions. I am also grateful to the Studienstiftung des Deutschen Volkes for their financial support.

\section{Engel structures and related distributions} \mlabel{s:defengel}

We first define all the distributions we will encounter in our constructions and the relations between them. This will lead to a proof of \propref{p:trivial intro}. To every Engel structure one can associate a foliation of rank $1$ and hypersurfaces transverse to the this foliation carry a contact structure together with a Legendrian line field. In \secref{ss:vert mod} we explain vertical modifications of the boundary, a construction which will be used frequently.
\subsection{Contact structures and even contact structures} \label{s:cont even cont}
Contact structures and even contact structures arise naturally on Engel manifolds. They will play an important role in all our constructions. Here we summarize the definitions and elementary properties. 
\begin{defn} \mlabel{d:cont str}
A {\em contact structure} $\CC$ on a $(2n-1)$--dimensional manifold $N$ is a smooth subbundle of $TN$ with corank $1$ such that for every local defining $1$--form $\alpha$ the rank of $d\alpha\eing{\CC}$ is maximal, or equivalently $\alpha\ww d\alpha^{n-1}\neq 0$. 
\end{defn}
If $n$ is even the sign of $\alpha\ww d\alpha^{n-1}$ is independent of the choice of $\alpha$. Thus a contact manifold of dimension $4n-1$ has a preferred orientation. In dimension $3$ the orientability of $M$ is the only obstruction to the existence of a contact structure by a result of J.~Martinet in \cite{martinet}.

The following two theorems will be used later. We sketch the proofs because we will apply them in concrete situations. First we explain Gray's theorem which states a remarkable stability property of contact structures.
\begin{thm}[Gray, \cite{gray}] \mlabel{t:gray}
Let $\CC_s, s\in[0,1]$ be a family of contact structures on $N$ which is constant outside a compact subset of $N$. Then there is an isotopy $\psi_s$ such that $\psi_{s*}\CC_0 = \CC_s$.
\end{thm}
\begin{proof}
We assume that $\CC_s$ is globally defined by a smooth family of $1$--forms $\alpha(s)$. The proof without this assumption is slightly more complicated, cf.~\cite{martinet}. The desired isotopy is the flow of a time--dependent vector field $Z(s)$. This is the unique vector field which is tangent to $\CC_s=\ker(\alpha(s))$ and satisfies
\begin{equation} \label{e:gray constr}
i_{Z(s)}d\alpha(s)=-\dot\alpha(s) \textrm{ on } \CC_s~.
\end{equation}
Because $d\alpha(s)$ is a non--degenerate $2$--form on $\CC_s$, such a vector field exists and is uniquely determined. One can show by a direct calculation that the flow of $Z(s)$ has the desired properties.
\end{proof}
An immediate consequence of Gray's theorem is a normal form for the contact structure on neighborhoods of curves tangent to the contact structure, i.e.~Legendrian curves. 
\begin{cor} \mlabel{c:tube leg curve}
If $\gamma$ is a closed Legendrian curve in a $3$--manifold with contact structure then we can choose coordinates $x,z,t$ on a tubular neighborhood of $\gamma$ such that $\gamma=\{x=z=0\}$ and the contact structure is defined by $dz-x\,dt$.
\end{cor}
In particular, every point of a contact manifold admits a neighborhood with coordinates $x,y,z$ such that the contact structure is defined by $\alpha=dz-x\,dy$.

Next we consider vector fields which preserve a given contact structure, so called contact vector fields. First recall that to each contact form $\alpha$ one can associate a unique contact vector field $R$, the {\em Reeb vector field}, such that $\al(R)\equiv 1$ and $i_{R}d\al\equiv 1$.
\begin{prop}  \mlabel{p:convectfunctions}
The map which assigns to each contact vector field $X$ the function $\alpha(X)$ is an isomorphism between the space of contact vector fields and $C^\infty(M)$.
\end{prop} 
\begin{proof}
Injectivity follows immediately from \defref{d:cont str}. Now let $f$ be a smooth function on $M$ and $R$ the Reeb vector field. Since $d\al\eing{\CC}$ is non--degenerate there is a unique vector field $Y$ tangent to $\CC$ such that $\left(i_Yd\al\right)\eing{\CC} = -df\eing{\CC}$. Then $fR+Y$ has the desired properties. This proves surjectivity.
\end{proof}
We will discuss several theorems from $3$--dimensional contact topology in \secref{s:contact top}. Now we state the definition and some elementary properties of even contact structures.
\begin{defn} \mlabel{d:evencontact}
Let $M$ be a $2n$--dimensional manifold and $\EE$ a distribution on $M$ of corank one. If for every local defining $1$--form $\alpha$ of $\EE$ the $2$--form $d\alpha$ has maximal rank on $\EE$, then $\EE$ is an even contact structure.
\end{defn}
Equivalently, $\EE$ is an even contact structure if for every local defining form $\alpha$, the $(2n-1)$--form $\alpha\ww d\alpha^{n-1}$ has no zeroes. In dimension $4$ this condition may be rephrased as $[\EE,\EE]=TM$. The definitions of contact structures and even contact structures are very similar. Still there are significant differences. One of them is the presence of a distinguished line field on a manifold with an even contact structure.
 
Since $\EE$ has rank $2n-1$, the rank of $d\alpha\eing{\EE}$ is $2n-2$. Hence $d\alpha\eing{\EE}$ has a kernel $\WW\subset\EE$ of dimension $1$ and it is easy to show that $\WW$ does not depend on the choice of the defining form $\alpha$. All flows which are tangent to $\WW$ preserve $\EE$. 
\begin{defn} \mlabel{d:charlinefield}
The line field $\WW$ is the {\em characteristic line field} of $\EE$. The foliation induced by this line field is called the {\em characteristic foliation}. A hypersurface is {\em transverse} if it is transverse to the characteristic foliation.  
\end{defn}
The following lemma is a simple consequence of the definition of an even contact structures and its characteristic foliation. 
\begin{lem} \mlabel{l:transverse contact structure}
Let $\EE$ be an even contact structure on $M$ and $\WW$ be the characteristic line foliation of $\EE$. If $N$ is a transverse hypersurface, then $TN\cap\EE$ is a contact structure on $N$. 

If $N'$ is another transversal such that two interior points $p\in N$ and $q\in N'$ lie on the same leaf $\WW_p$ of the characteristic foliation, then the map obtained by following nearby leaves preserves the induced contact structures on neighborhoods of $p$ and $p'$
\end{lem}
As we have explained above a contact structure on a manifold of dimension $2n-1$ induces an orientation of this manifold if $n$ is even. In this situation an orientation of the characteristic foliation together with the contact orientation on a transversal define an orientation of $M$.
\begin{prop} \mlabel{p:even contact orientations}
Let $\EE$ be an even contact structure on a $4n$--manifold $M$. Then an orientation of $M$ induces an orientation of the characteristic line field $\WW$ and vice versa.
\end{prop}
  
\subsection{Engel structures -- Definition and first examples} \mlabel{s:engeldef}
 
If $\DD$ is a distribution we denote the sheaf of all vectors which can be obtained as commutators of pairs of local sections of $\DD$ by $[\DD,\DD]$.   
\begin{defn} \mlabel{d:Engel}
An Engel structure is a distribution $\DD$ of rank two on a manifold $M$ of dimension four with the following properties.
\begin{itemize}
\item[(i)] $\EE=[\DD,\DD]\subset TM$ is a subbundle of rank $3$.
\item[(ii)] $TM=[\EE,\EE]$.
\end{itemize}
\end{defn}
The second condition in \defref{d:Engel} implies that $\EE=[\DD,\DD]$ is an even contact structure. To $\EE$ we associate the characteristic foliation $\WW$ of $\EE$. By the definition of $\WW$ we have $\WW\subset\EE$. 
\begin{lem} \mlabel{l:EL}
 If $\EE=[\DD,\DD]$ is an even contact structure induced by Engel structure $\DD$, then $\WW\subset\DD$.
\end{lem}
\begin{proof}
Assume that $\WW_p$ is not contained in $\DD_p$ at $p\in M$. Let $\alpha$ be a local defining form for $\EE$ on a neighborhood of $p$. If $X,Y$ are two linearly independent local sections of $\DD$ around $p$, then $d\alpha(X_p,Y_p)\neq 0$ by the assumption on $\WW$. On the other hand $d\alpha(X,Y) = -\alpha([X,Y])=0$ since $[X,Y]$ is a local section of $\EE=[\DD,\DD]$. Thus the assumption $\WW_p\not\subset\DD_p$ leads to a contradiction.
\end{proof}
Together with \lemref{l:transverse contact structure} this allows us to give the following interpretation of the condition $\EE=[\DD,\DD]$. As one moves along the leaves of $\WW$, the plane field $\DD$ rotates in $\EE$ around the characteristic foliation without stopping. 

We have shown that an Engel structure $\DD$ on $M$ induces a flag of distributions
\begin{equation}  \label{e:flag of distr}
\WW\subset\DD\subset\EE\subset TM~.
\end{equation}
Each of these distributions has corank one in the distribution containing it. By \propref{p:even contact orientations}, an orientation of the characteristic foliation of an Engel structure induces an orientation of the underlying manifold and vice versa. In addition, $\EE=[\DD,\DD]$ is oriented by the following proposition.
\begin{prop} \mlabel{p:E is oriented}
If $\DD$ is an Engel structure, then the even contact structure $\EE=[\DD,\DD]$ has a distinguished orientation. 
\end{prop}
\begin{proof}
Let $X,Y$ be local sections of $\DD$ around $p\in M$ such that $X(p)$ and $Y(p)$ are linearly independent. Then $X(p), Y(p), [X,Y](p)$ is an orientation of $\EE(p)$ which is independent of the choice of $X,Y$. 
\end{proof}
This leads to the following orientation conventions.
\begin{itemize}
\item[(i)] If $\WW$ is oriented we orient $M$ by \propref{p:even contact orientations}.
\item[(ii)] $\EE=[\DD,\DD]$ is oriented according to \propref{p:E is oriented} .
\item[(iii)] Transverse hypersurfaces are oriented by the induced contact structure. In particular transverse boundaries are oriented by this convention.
\item[(iv)] If $\WW$ is oriented, then we orient the contact structure on a transverse hypersurface such that the orientation of the contact structure followed by the orientation of $\WW$ is the orientation of the even contact structure. 
\end{itemize}
Let $M$ be an Engel manifold with oriented characteristic foliation and transverse boundary. We write $\partial_+M$ for those connected components of $\partial M$ where the characteristic foliation points out of $M$. The remaining  boundary components are $\partial_-M=\partial M\setminus \partial_+M$. 

It is clear from \eqref{e:flag of distr} that the existence of an Engel structure has strong implications for the topology of the underlying manifold. A proof of the following result can be found in \cite{montgomery2}. It was known earlier by V.~Gershkovic.
\begin{prop} \mlabel{p:trivial TM}
Let $\DD$ be an oriented Engel structure on an oriented $4$--manifold $M$. Then the tangent bundle of $M$ is trivial.
\end{prop}
\begin{proof}
Consider the flag $\WW\subset \DD \subset \EE \subset TM$ of subbundles of $TM$. According to our orientation conventions each of the quotient bundles $\DD/\WW$, $\EE/\DD$ and $TM/\EE$ is oriented. Thus
\begin{equation} \label{e:TMsumme}
TM = \WW\oplus\frac{\DD}{\WW}\oplus\frac{\EE}{\DD}\oplus\frac{TM}{\EE}
\end{equation}
is isomorphic to the sum of four trivial line bundles. 
\end{proof}
Under the assumptions of \propref{p:trivial TM} we can choose oriented orthogonal complements in \eqref{e:flag of distr} after having fixed an auxiliary Riemannian metric and we obtain a trivialization of $TM$. Since we only chose the Riemannian metric, this framing is well defined up to homotopy. We will refer to such framings as {\em Engel framings}.

A very simple example of an Engel structure is the standard Engel structure in the following Darboux theorem for Engel structures. 
\begin{thm}[Engel, \cite{engel}]
Every point of an Engel manifold has a neighborhood with coordinates $w,x,y,z$ such that the Engel structure is given by
\begin{equation} \label{e:standard R4}
\ker(dz-x\,dy)\cap\ker(dx-w\,dy)~.
\end{equation}
\end{thm}
Here, the characteristic foliation is spanned by $\partial_w$ and the even contact structure is defined by $dz-x\,dy$. We now discuss more interesting examples. Apart from the construction we present later, two other constructions of closed Engel manifolds can be found in the literature.
\begin{ex} \mlabel{e:prolong}
The starting point of the {\em prolongation} construction is a contact structure $\CC$ on a $3$--manifold $N$. The projectivization $\PP\CC$ of $\CC$ is a circle bundle over $N$. The projection $\pr : \PP\CC\rightarrow N$ of this fibration maps each Legendrian line to its base point in $N$. One can show that
$$
\DD_\CC = \left\{ v\in T_{[l]}\PP\CC\: \big|\: \pr_*(v)\in [l] \right\}
$$
is an Engel structure whose characteristic foliation is tangent to the fibers of $\pr$. The associated even contact structure is $\pr_*^{-1}\CC$. Similarly, one obtains {\em orientable} Engel structures on the space of {\em oriented} Legendrian lines.
\end{ex}
The second previously known construction of closed Engel manifolds is due to H.--J.~Geiges. It yields Engel structures on parallelizable mapping tori of compact $3$--manifolds, cf.~\cite{geiges}. 

\subsection{Transverse hypersurfaces in Engel manifolds} \mlabel{s:int line field} 
According to \lemref{l:transverse contact structure} transverse hypersurfaces carry a natural contact structure. In this section we discuss a Legendrian line field induced by the Engel structure on a transverse hypersurface. This line field can be used to obtain a normal form for the Engel structure near a transverse hypersurface $N$.
\begin{defn}
The {\em intersection line field} of $\DD$ on $N$ is the Legendrian line field $TN\cap\DD$.
\end{defn}
If the Engel structure $\DD$ is oriented, then we orient the intersection line field such that the orientation of $\WW$ followed by the orientation of the intersection line field is the orientation of $\DD$.

The following theorem, as well as \propref{p:Engel diffeos}, can be found in \cite{montgomery} but according to that article they were known before. Both deal with Engel structures obtained by prolongation.

If one applies the prolongation construction to the contact structure $\CC$ on a transverse hypersurface $N$ in an Engel manifold, then one obtains the manifold $\PP\CC$ with its canonical Engel structure. Let $\sigma$ be the section of $\PP\CC$ which assigns to each point $p$ of $N$ the Legendrian line $TN\cap\DD(p)$.
\begin{thm} \mlabel{t:normal forms for transversals}
Any sufficiently small tubular neighborhood of $N$ in $M$ is canonically diffeomorphic as an Engel manifold to a tubular neighborhood of $\sigma$ in $\PP\CC$.
\end{thm}
This means that the intersection line field and the contact structure $\CC$ on a transverse hypersurface $N$ suffice to recover the Engel structure on a small neighborhood of $N$. 

Another property of prolonged Engel structures is that they admit many diffeomorphisms which preserve the Engel structure, cf.~\cite{montgomery}. 
\begin{prop} \mlabel{p:Engel diffeos}
Let $N_1,N_2$ be $3$--manifolds with contact structures $\CC_1,\CC_2$ and let $\varphi : N_1 \rightarrow N_2$ be a contact diffeomorphism. The diffeomorphism between $\PP\CC_1$ and $\PP\CC_2$ which is induced by the action of $\varphi$ on $\CC_1$ preserves the Engel structures. 

Every diffeomorphism $\PP\CC_1 \rightarrow \PP\CC_2$ with this property is of this form.
\end{prop}
We will mainly be concerned with the case that $N$ is the transverse boundary of a compact Engel manifold. In the following section we will explain a modification of the Engel structure on a collar of the boundary which allows us to vary the intersection line field within its homotopy class. Therefore it suffices to control the homotopy class of the intersection line field and not the induced foliation. 

In the remaining part of this section we define simple invariants which determine the intersection line field on transverse boundaries of Engel manifolds. For this we assume that the characteristic foliation and the Engel structure are oriented. This induces orientations of the contact structure and the intersection line field on $N$. 

If we fix an oriented framing $C_1,C_2$ of the contact structure, we can assign an element of $H^1(N;\Z)$ to a nowhere vanishing Legendrian line field $X$ as follows. There are unique functions $f_1,f_2$ such that $X=f_1C_1+f_2C_2$. The element of $H^1(N;\Z)$ corresponding to $X$ assigns the winding number of  
$$
(f_1\circ\gamma, f_2\circ\gamma) :  S^1\lra \R^2\setminus\{0\}
$$
around $0$ to the homology class represented by $\gamma$. This depends only on the homotopy class of $X$ as nowhere vanishing Legendrian line field.

In general we will not use a distinguished global framing but we will compare the homotopy classes of nowhere vanishing Legendrian vector fields using rotation numbers along Legendrian curves. Let $X$ be a nowhere vanishing Legendrian vector field on $N$ and let $X,Y$ be an oriented framing of the contact structure on $N$. If $\gamma$ is an immersed, closed Legendrian curve, then $\dot\gamma$ is a Legendrian vector field along $\gamma$ and there are unique functions $f_1,f_2$ such that $\dot\gamma=f_1X +f_2Y$.
\begin{defn} \mlabel{d:rotzahl}
The winding number of the map $(f_1\circ\gamma,f_2\circ\gamma) : S^1\lra \R^2\setminus\{0\}$ around zero is the {\em rotation number} of $X$ along $\gamma$. 
\end{defn}
The sign of the rotation number changes if the orientation of $\gamma$, or if the orientation of the contact structure is reversed. In particular, its sign changes if we change the orientation of the characteristic foliation of $\DD$ but it is independent of the orientation of $\DD$.  

In order to determine the homotopy class of the intersection line field, it suffices to know the rotation number along sufficiently many Legendrian curves.

\subsection{Vertical modifications of transverse boundaries} \mlabel{ss:vert mod}

In this section we consider a manifold $M$ with oriented Engel structure $\DD$ with oriented characteristic foliation such that $\partial M$ is transverse and compact. Vertical modifications of the boundary allow us to change the intersection line field on $\partial M$. In this construction we modify the Engel structure near the boundary such that the even contact structure remains unchanged. Assume that $\LL$ is a Legendrian line field which is homotopic to the intersection line field $\LL_\DD$ of $\DD$. Then we can modify $\DD$ such that the intersection line field of the new Engel structure is $\LL$.

This can be generalized to the case of a transverse hypersurface $N$, but then the possible changes of the intersection line fields depend on the behavior of the Engel structure along the leaves of the characteristic foliation through $N$. This problem will arise only in \thmref{t:M+N zweite vers.} at the end of this article. 

Let $X_0$ be a Legendrian vector field which spans and orients the intersection line field on $\partial M$. Assume that the Legendrian vector field $X_1$ on $\partial M$ is homotopic to $X_0$ through a family $X_s, s\in[0,1]$ of nowhere vanishing Legendrian vector fields. Our aim is to construct an Engel structure on $M$ such that the intersection line field on $\partial M$ is spanned and oriented by $X_1$ without changing the even contact structure. We treat the boundary components $\partial_+M$ where $\WW$ points out of $M$. The components $\partial_-M=\partial M\setminus\partial_+M$ can be treated similarly.

Choose a collar $U\simeq\partial_+M\x (-1,0]$ of $\partial_+M$ such that the second factor with its standard orientation corresponds to the oriented characteristic foliation of $\DD$. We write $w$ for the coordinate on $(-1,0]$.

Furthermore let $Y_0$ be a vector field such that $X_0,Y_0$ is an oriented trivialization of the contact structure on $\partial_+M$. The even contact structure on $U$ is spanned and oriented by $X_0,Y_0,\partial_w$. There is a unique smooth function $f : U \simeq \partial_+M\x(-1,0]\lra \R$ which vanishes along $\partial_+M$ such that $\partial_w$ and
\begin{equation} \label{e:X}
X(p,w)=\cos(f(p,w)) X_0(p) + \sin(f(p,w)) Y_0(p)
\end{equation}
span and orient $\DD(p,w)$. Because $\DD$ is an Engel structure $X,[\partial_w,X]$ and $\partial_w$ are linearly independent everywhere. This implies that $f$ is strictly monotone along the leaves of $\WW$. By our orientation convention for the contact structure on transverse boundaries, $f$ is strictly increasing. 

There is a unique family of real functions $g_s$ such that $\cos(g_s(p))X_0+\sin(g_s(p))Y_0$ is a positive multiple of $X_s$ and $g_0=0$. Since $\partial_+M$ is compact we can choose $k\in\Z$ such that $g_1+2\pi k$ is at least $2\pi$. Now extend $f$ from $\partial_+M\x(-1,0]$ to a smooth function $\widehat{f}$ on $\partial_+M\x(-1,1]$ such that $\widehat{f}$ is strictly increasing and coincides with $g_1+2\pi k$ along $\partial_+M\x\{1\}$. We extend $\DD$ from $M$ to a distribution on $M\cup_\partial\partial_+M\x[0,1]$ by the span of the vector field $\partial_w$ and
$$
X(p,w)=\cos(\widehat{f}(p,w))X_0(p)+\sin(\widehat{f}(p,w))Y_0(p)~.
$$
The fact that $\partial_w\widehat{f}>0$ implies that $X$ and $[\partial_w,X]$ are linearly independent. Moreover, $\partial_w$ together with $X,[\partial_w,X]$, respectively the horizontal lifts of $X_0,Y_0$, span the same even contact structure. Thus the extended distribution is an Engel structure whose intersection line field on $\partial_+M\x\{1\}$ is spanned by $X_1$. Using a flow tangent to the characteristic foliation one can identify $M$ and $M\cup_\partial \partial_+M\x[0,1]$ such that the even contact structures are preserved.
\begin{defn}
The modifications of an Engel manifold with boundary described above will be called a {\em vertical modifications of the boundary}.
\end{defn}
Let us summarize the discussion in the following theorem. 
\begin{thm}   \mlabel{t:inters fol}
Let $(M,\DD)$ be an oriented Engel manifold such that $\partial M$ is transverse and compact. Assume that the Legendrian line field $\LL$ on $\partial_+M$ is homotopic to the intersection line field $\LL_\DD$ of $\DD$. 

Then $\LL$ is the intersection line field of a suitable vertical modification of $\partial_+M$.
\end{thm}
There is an analogous construction in the context of contact structures. Let $N$ be a compact $3$--manifold carrying a contact structure $\CC$ and a contact vector field $V$ transverse to the boundary. For this, let $X$ be a nowhere vanishing Legendrian vector field such that $[V,X]$ and $X$ are linearly independent. Assume that $X_1$ is a Legendrian vector field along $\partial N$ which is homotopic to the restriction of $X$ to $\partial N$ through nowhere vanishing Legendrian vector fields. Then we can modify $X$ near the boundary of $N$ such that the resulting Legendrian vector field $X'$ coincides with $X_1$ along $\partial  N$ and $X',[V,X']$ are linearly independent. 
\begin{rem} \mlabel{r:suitable vert mod 1}
In our applications the Legendrian vector field $X_1$ is given only on a compact embedded submanifold $U\subset\partial_+M$. We then pick a positive function $g_1+2\pi k$ as above and extend it to a positive function $\widetilde{g}$ which equals $2\pi$ outside of a small neighborhood of $U$. If $U$ is not connected we can choose different values for $k$ on each connected component of $U$. 
\end{rem}
%%%%%%%%%%%%%%%%%%%%%%%%%%%%%%%%%%%%%%%%%%%%%%%%%%%%%%%%%%%%%%%%%%%%%%%%%%
\section{Contact topology} \mlabel{s:contact top}

In this section we summarize facts from the theory of contact structures. \secref{s:stabilization} deals with Legendrian curves in contact manifolds. We discuss in particular stabilization of Legendrian curves. This will play an important role in the construction of attaching maps for round handles of index $1$ carrying model Engel structures. 

In \secref{s:contact facts} we collect those facts about convex surfaces in contact manifolds which we will use in \secref{s:bypasses} for the construction of bypasses in overtwisted contact manifolds as well as in \secref{s:modelle index 2} and \secref{ss:model 3-handles} for the construction of model Engel structures on round handles of index $2$ and $3$. 

Most of the material presented in this section (with the exception of the construction of bypasses in overtwisted contact manifolds in \secref{s:bypasses}) can be found in \cite{sympgeo, etnyrehonda1, etnyrehonda2, giroux, honda}.

\subsection{Properties of Legendrian curves} \mlabel{s:stabilization}

Throughout this section we consider an oriented contact structure $\CC$ on a $3$--manifold $N$. A curve is {\em Legendrian} if it is tangent to the contact structure.
\begin{defn} \mlabel{d:stand framing}
A framing $(S,T)$ of an embedded closed Legendrian curve $\gamma$ is an {\em oriented contact framing} if
\begin{itemize}
\item[(i)] $S$ is tangent to the contact structure and $\dot\gamma,S$ represents the orientation of the contact structure and
\item[(ii)] $T$ is transverse to the contact structure and $\dot\gamma,S,T$ represents the orientation induced by the contact structure.
\end{itemize}
\end{defn} 
Every Legendrian curve has a contact framing which is well defined up to homotopy. %If two framings $(S,T)$ and $(S',T')$ of $\gamma$ are homotopic we write $(S,T)\sim(S',T')$. 
The homotopy class of framings of $\gamma$ represented by contact framings is denoted by $\fr(\gamma)$. On the set of framings of $\gamma$ there is a $\Z$--action defined by 
\begin{align} \label{e:action on framings}
\begin{split}
\big(m\cdot(S,T)\big)(\gamma(t)) & = \big(\cos(mt)S(\gamma(t)) + \sin(mt)T(\gamma(t)), \\
                                 & \qquad -\sin(mt)S(\gamma(t)) + \cos(mt)T(\gamma(t))\big)~.
\end{split}
\end{align}
This action is transitive on the set of homotopy classes of framings such that $\dot\gamma$ followed by the framing induces the contact orientation. When we reverse the orientation of $\gamma$ the coorientation of $\gamma$ changes. Therefore the $\Z$--action on the framings does not depend on the orientation of $\gamma$. Contact framings are related to the Thurston--Bennequin invariant of null--homologous Legendrian knots. 
\begin{defn}  \mlabel{d:tb}
Let $\gamma$ be a null--homologous Legendrian curve in $N$. Fix a relative homology class $[\Sigma]\in H_2(N,\gamma;\Z)$ which is represented by an oriented surface $\Sigma$ such that $\partial \Sigma=\gamma$ and $\gamma$ is oriented as the boundary of $\Sigma$. A new curve $\gamma'$ is obtained by pushing $\gamma$ slightly along a vector field which is transverse to $\CC$. The {\em Thurston--Bennequin invariant} $\tb(\gamma,[\Sigma])$ is the homological intersection number of $\gamma'$ with $\Sigma$. 
\end{defn}
In \defref{d:tb} the surface $\Sigma$ singles out those framings of $\gamma=\partial N$ whose first component is tangent to $\Sigma$. We write $\fr_\Sigma(\gamma)$ for the corresponding homotopy class of framings. The Thurston--Bennequin invariant measures the difference between $\fr_\Sigma(\gamma)$ and $\fr(\gamma)$ of $\gamma$
\begin{equation} \label{e:seifert tb}
\tb(\gamma,[\Sigma])\cdot \fr_\Sigma(\gamma) = \fr(\gamma)~.
\end{equation}
A simple application of Gray's theorem (\thmref{t:gray}) shows that the contact framing of a Legendrian curve determines the isotopy type of the contact structure on a neighborhood of the Legendrian curve.
\begin{lem} \mlabel{l:tub nbhd of leg curves}
Let $\CC_0,\CC_1$ be two contact structures on $N$ which induce the same orientation of $M$ and $\gamma$ a curve tangent to both $\CC_0$ and $\CC_1$. If the contact framings of $\gamma$ with respect to $\CC_0$ and $\CC_1$ are homotopic, then there is an isotopy $\psi_s, s=[0,1],$ of $N$ relative $\gamma$ such that $\psi_{1*}\CC_0=\CC_1$ on a tubular neighborhood of $\gamma$.

If $\CC_0$ and $\CC_1$ are oriented, then one can choose $\psi_s$ such that $\psi_1$ preserves also the orientations of the contact structures.
\end{lem}
Now we discuss the other classical invariant of null--homologous Legendrian curves.
\begin{defn}  \mlabel{d:rot zahl klassisch}
Let $\Sigma$ be a connected orientable surface with $\partial\Sigma=\gamma$. Fix an oriented trivialization $X,Y$ of $\CC$ along $\Sigma$. There are unique functions $f_x,f_y$ such that $\dot\gamma(t)=f_X(t)X + f_Y(t)Y$. The winding number of
\begin{align*}
S^1 & \lra \R^2\setminus\{0\} \\
  t & \lmt (f_Y(t),f_Y(t))
\end{align*}
around $0$ is the {\em rotation number} $\rot(\gamma,[\Sigma])$.
\end{defn}

\begin{rem} \mlabel{r:rotzahlvergleich}
In \defref{d:rotzahl} we defined the rotation number for Legendrian curves on transverse hypersurfaces in oriented manifolds with oriented Engel structures. In \defref{d:rot zahl klassisch} we fix an oriented trivialization of the oriented contact structure on a Seifert surface $\Sigma$ of the Legendrian knot $\partial\Sigma=\gamma$ and compare $\dot\gamma$ with this trivialization. 

Now assume that $\CC$ is the oriented contact structure on a transverse hypersurface of an oriented manifold with oriented Engel structure. In this case we can use the intersection line field as the first component of the trivialization of $\CC$ over $\Sigma$. Thus in this situation the two rotation numbers in \defref{d:rotzahl} and \defref{d:rot zahl klassisch} are equivalent.
\end{rem}
Using contact framings and rotation numbers one can distinguish isotopy classes of Legendrian curves. This isotopy type can be changed using the following construction. 

\begin{figure}[htb]
\begin{center}
\includegraphics[width=0.96\textwidth]{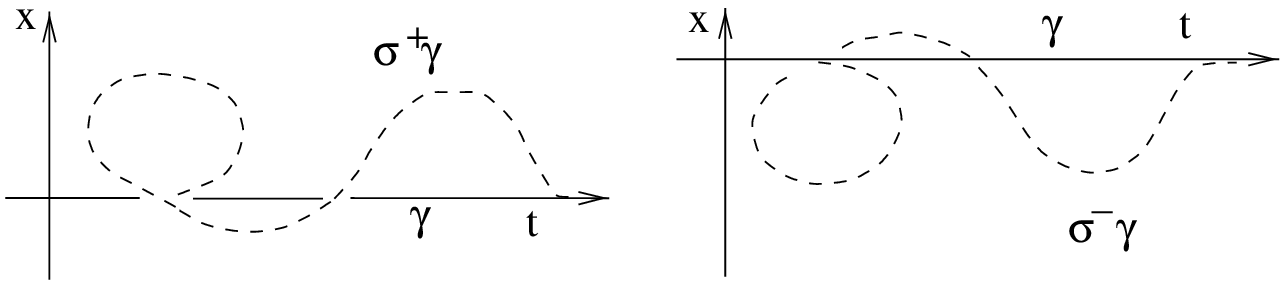}
\caption{\label{b:stab}}
\end{center}
\end{figure}

Let $\gamma$ be a closed Legendrian curve in a contact manifold. According to \corref{c:tube leg curve} there are coordinates $x,z,t$ on a neighborhood of $\gamma$ such that the contact structure is defined by $dz-x\,dt$ and $\gamma=\{x=0, z=0\}$. The $t$--axis in both parts of \figref{b:stab} represents the projection of a neighborhood of a segment of $\gamma$ to the $x,t$--plane and $\partial_z$ is pointing outwards. Let $\sigma$ be another smooth curve in the $t,x$--plane that coincides with the $t$--axis near its endpoints. Since the contact structure near $\gamma$ is defined by $dz-x\,dt$ the segment $\sigma$ lifts to a Legendrian curve whose endpoints lie on $\gamma$ if and only if the integral of the $1$--form $x\,dt$ along $\sigma$ is zero. Then the lifts of curves $\sigma^+\gamma$, respectively $\sigma^-\gamma$ as in \figref{b:stab}, yield new closed Legendrian curves which are close to $\gamma$. They represent the positive and the stabilization of $\gamma$. 

Up to isotopy through Legendrian curves the stabilized curves do not depend on the choices made above. In particular, positive and negative stabilization commute up to Legendrian isotopy. The effect of stabilization on the rotation number with respect to a nowhere vanishing Legendrian vector field $X$ is given by
\begin{align} \label{e:rot stab}
\begin{split}
\rot(\sigma^+\gamma,X) & = \rot(\gamma,X)+1 \\
\rot(\sigma^-\gamma,X) & = \rot(\gamma,X)-1~.
\end{split}
\end{align}
Hence $\sigma^+\gamma,\sigma^-\gamma$ and $\gamma$ are pairwise non--isotopic as Legendrian curves. From the construction of $\sigma^\pm\gamma$ it follows that
\begin{equation} \label{e:tb twist}
\fr(\sigma^+\gamma) = 1\cdot\left(\psi_{1*}\fr_\CC(\gamma)\right) = \fr(\sigma^-\gamma)~. 
\end{equation}
From \eqref{e:seifert tb} it follows that this is coherent with the well known formula $\tb(\sigma^\pm\gamma)=\tb(\gamma)-1$ for null--homologous Legendrian curves.

%%%%%%%%%%%%%%%%%%%%%%%%%%%%%%
\subsection{Facts from the theory of convex surfaces} \mlabel{s:contact facts}

Let $N$ be a $3$--manifold with an oriented contact structure $\CC$. Consider a properly embedded oriented surface $\Sigma$. The contact structure induces the {\em singular foliation} $\FF=\CC\cap T\Sigma$ on $\Sigma$. Usually this is called the characteristic foliation of $\Sigma$. We avoid this terminology since there is already a characteristic (non--singular) foliation in the context of Engel structures. 

The singularities of $\FF$ are those points $p\in\Sigma$ where $\CC_p=T_p\Sigma$. From the orientations of $\Sigma$ and $\CC$ we obtain an orientation of $\FF$ away from the singularities using the following convention: For each non--singular point $p$ of $\FF$ we choose 
$$
v\in\FF_p, v_\Sigma\in T_p\Sigma\setminus\FF_p \textrm{ and } v_\CC\in\CC_p\setminus\FF_p
$$
such that $(v,v_\CC,v_\Sigma)$ is the contact orientation, $(v,v_\Sigma)$ orients $\Sigma$ and $(v,v_\CC)$ orients $\CC$. Then $v$ represents the orientation of $\FF_p$. 

Generically, singular points are non--degenerate. We say that a singular point is {\em elliptic} if its index is $+1$ and {\em hyperbolic} if the index is $-1$. When the orientation of $\CC$ and the orientation of the surface coincide at a singular point of $\FF$, we say that this singularity is {\em positive}, otherwise it is {\em negative}. If we orient $\FF$ according to our conventions, positive elliptic points are sources and negative elliptic points are sinks.
\begin{defn}
$\Sigma$ is called {\em convex} if there is a contact vector field transverse to $\Sigma$. 
\end{defn}
Giroux studied closed convex surfaces in \cite{giroux}. If $\Sigma$ has boundary we will usually assume that $\partial\Sigma$ is Legendrian. In particular he showed that a closed embedded surface is generically convex (with respect to the $C^\infty$--topology).  The analogous statement is in general not true when $\Sigma$ has boundary (even when $\partial \Sigma$ is Legendrian). For each boundary component $\gamma\subset\partial\Sigma$ let $t(\gamma,\fr_\Sigma)$ be the unique integer such that 
$$
t(\gamma,\fr_\Sigma)\cdot\fr_\Sigma(\gamma)=\fr(\gamma)~. 
$$
If $\gamma$ is a Legendrian knot and $\Sigma$ is a Seifert surface for $\gamma$, then $t(\gamma,\fr_\Sigma)$ is the Thurston--Bennequin invariant by \eqref{e:seifert tb}. 
\begin{prop}[Honda, \cite{honda}] \mlabel{p:convexity with boundary}
Let $\Sigma$ be a compact, oriented, properly embedded surface with Legendrian boundary, and assume $t(\gamma,\fr_{\Sigma}) \le 0$ for all boundary components $\gamma$ of $\Sigma$. There exists a $C^0$--small perturbation near the boundary (fixing $\partial\Sigma$) which puts an annular neighborhood $A$ of $\partial \Sigma$ into a standard form, and a subsequent $C^\infty$--small perturbation keeping $A$ fixed which makes $\Sigma$ convex. 

Moreover, if $V$ is a contact vector field defined on a neighborhood of $A$ and transverse to $A\subset \Sigma$, then $V$ can be extended to a contact vector field transverse to all of $\Sigma$. 
\end{prop}
The singular foliation is enough to determine the contact structure on a small neighborhood of a convex surface. 
\begin{thm}[Giroux, \cite{giroux}] \mlabel{t:giroux unique}
Let $\Sigma$ be a compact orientable convex surface with Legendrian boundary. Two $\R$--invariant contact structures on $\Sigma\x\R$ which induce the same orientation and the same singular foliation on $\Sigma\x\{0\}$ are isotopic. They are conjugate by a diffeomorphism $\varphi\x\id$, and $\varphi$ is isotopic to the identity through diffeomorphisms of $\Sigma$ that preserve the singular foliation.
\end{thm}
For a convex surface $\Sigma$ with Legendrian boundary we fix a contact vector field $V$ transverse to $\Sigma$. 
\begin{defn} \mlabel{d:dividing set}
The {\em dividing set} $\Gamma_\Sigma$ of $\Sigma$ is the set of points where $V$ is tangent to the contact structure.
\end{defn}
The dividing set contains much information about the contact structure near $\Sigma$. Giroux showed in \cite{giroux} that $\Gamma_\Sigma$ is a submanifold of $\Sigma$ which is transverse to the singular foliation. Its isotopy class depends only on $\Sigma$ and $\CC$ but not on $V$. 
\begin{defn}
Let $\FF$ be a singular foliation on $\Sigma$ such that $\partial\Sigma$ is tangent to $\FF$. A collection $\Gamma\subset\Sigma$ of closed curves and arcs with end points on $\partial\Sigma$ is said to {\em divide} $\FF$ if on each connected component of the closure of $\Sigma\setminus\Gamma$ there is a smooth volume form $\omega$ and a vector field $X$ tangent to $\FF$ such that 
\begin{itemize}
\item[(i)] the divergence of $X$ with respect to $\omega$ is positive everywhere
\item[(ii)] $X$ points outwards along those parts of the boundary which correspond to $\Gamma$. 
\end{itemize} 
\end{defn}
\begin{thm}[Giroux, \cite{giroux}] \mlabel{t:divide impera}
If $\Sigma$ is a convex surface with Legendrian boundary, then $\Gamma_\Sigma$ divides the singular foliation on $\Sigma$.

Conversely, let $\FF$ be a singular foliation on a compact oriented surface $\Sigma$ such that $\partial\Sigma$ is tangent to $\FF$ and $\Gamma$ divides $\FF$. Then there is a positive $\R$--invariant contact structure on $\Sigma\x\R$ such that $\Sigma\x\{0\}$ is convex, the induced singular foliation on $\Sigma\x\{0\}$ is precisely $\FF$ and such that $\Gamma$ consists of those points where the contact structure is tangent to the second factor in $\Sigma\x\R$.
\end{thm}
This theorem implies in particular that the dividing set of a convex surface with Legendrian boundary is not empty. Assume that $\CC$ is cooriented  by $\alpha$ and that the orientation $\Sigma$ followed by $V$ is the contact orientation. Then the dividing set $\Gamma$ separates the region $\Sigma_+$ where $\alpha(V)$ is positive from the region $\Sigma_-$ where $\alpha(V)$ is negative. If $\Sigma$ is closed, then the Euler characteristic $\chi(\CC)$ of $\CC$ viewed as oriented bundle satisfies
\begin{equation} \label{e:eulerchar}
\chi(\CC) = \chi(\Sigma_+)-\chi(\Sigma_-)~.
\end{equation}
If $\Sigma$ is a Seifert surface for a Legendrian knot $\gamma$, then one can determine the Thurston--Bennequin invariant and the rotation number of $\gamma$ using the following formulas from \cite{kanda2, honda}
\begin{align} \label{e:kanda rot}
\begin{split}
\tb(\gamma) & = -\frac{1}{2}\left|\Gamma\cap\gamma\right| \\
\rot(\gamma)& = \chi(\Sigma_+)-\chi(\Sigma_-)~.
\end{split}
\end{align} 
Next we consider deformations of the singular foliation. Let $\Sigma$ be a convex surface with Legendrian boundary and fix a transverse contact vector field $V$. 
\begin{defn}
An isotopy $\Phi_s$ of a surface $\Sigma$ is called {\em admissible} if $\Phi_s(\Sigma)$ is transverse to $V$ for all $s$. 
\end{defn}
The following theorem is a generalization of the Giroux flexibility theorem in \cite{giroux} where $\Sigma$ is assumed to be closed.
\begin{thm}[Giroux, Honda,\cite{giroux, honda}] \mlabel{t:giroux flex}
Let $\FF_0$ be the singular foliation on $\Sigma$ induced by the contact structure and $\FF_1$ a singular foliation which is divided by $\Gamma_\Sigma$. Then there is an admissible isotopy $\Phi_s,s\in[0,1]$, of $\Sigma$ such that $\FF_1$ is the singular foliation on $\Phi_1(\Sigma)$.
\end{thm}
\begin{ex} 
In this example we want to fix some terminology. Consider the $\R$--invariant contact form $\alpha= \cos(\varphi)\,dt +\sin(\varphi)\,dx$ on $T^2\x\R$ where $x$ is the coordinate on the $\R$--factor while $\varphi$ and $t$ correspond to $T^2=S^1\x S^1$. We say that the singular foliation $\FF$ on $T^2\x\{0\}$ is in {\em standard form}. The singularities of the singular foliation form two circles $\{\varphi=\pi/2\}\cup\{\varphi=3\pi/2\}$. These are called the {\em Legendrian divides}. The dividing set of $T^2\x\{0\}$ has two connected components
$$
\Gamma = \{ \varphi = 0\} \cup \{\varphi =\pi\}~.
$$
The curves tangent to $\partial_\varphi$ form the {\em Legendrian ruling}. By \thmref{t:giroux flex}, the Legendrian ruling can be changed into any foliation by straight lines as long as the straight lines remain transverse to the dividing set.
\end{ex}
Using \thmref{t:giroux flex} we can isotope $\Sigma$ such that a submanifold with certain properties becomes Legendrian. 
\begin{defn}
A union $C$ of disjoint properly embedded arcs and closed curves on $\Sigma$ is called {\em non--isolating} if
\begin{itemize}
\item[(i)] $C$ is transverse to $\Gamma$ and every arc begins and ends on $\Gamma$
\item[(ii)] the boundary of each component of $\Sigma\setminus(\Gamma\cup C)$ contains a segment of $\Gamma$. 
\end{itemize}
\end{defn}
\begin{thm}[Kanda, Honda, \cite{kanda, honda}]  \mlabel{t:LeRP}
Consider a collection $C$ of properly embedded closed curves and arcs on a convex surface $\Sigma$ with Legendrian boundary. If $C$ is non--isolating there exists a $C^0$--small admissible isotopy $\Phi_s, s\in[0,1]$ such that 
\begin{itemize}
\item[(i)] $\Phi_0=\mathrm{id}$ and $\Phi_1(\Gamma_\Sigma)=\Gamma_{\Phi_1(\Sigma)}$
\item[(ii)] $\Phi_1(C)$ is Legendrian.
\end{itemize}
\end{thm}
We finish this section with an important dichotomy of contact structures on $3$--manifolds. 
\begin{defn}
An embedded disc with Legendrian boundary is {\em overtwisted disc} if all singularities on the boundary have the same sign. A contact structure is called {\em overtwisted} if it there is an overtwisted disc. A contact structure is {\em tight} if it is not overtwisted.
\end{defn}
The following theorem shows that tight contact structures are more interesting than overtwisted ones in many aspects.
\begin{thm}[Eliashberg, \cite{eliashot}] \mlabel{t:eliashbergs wahnsinn}
If two overtwisted contact structures on a closed manifold are homotopic as plane fields and induce the same orientation, then they are isotopic.
\end{thm}
For our purposes however, the flexibility of overtwisted contact structures will be very useful. Regions where a given contact structure is tight can be found using the following theorem.
\begin{thm}[Colin, \cite{colin}] \mlabel{t:giroux tight}
If $\Sigma\neq S^2$ is a convex surface (closed or compact with Legendrian boundary) in a contact manifold $(M,\CC)$, then $\Sigma$ has a tight neighborhood if and only if the dividing set of $\Sigma$ has no homotopically trivial closed curves. If $\Sigma=S^2$, then $\Sigma$ has a tight neighborhood if and only if $\Gamma_\Sigma$ is connected.
\end{thm}

\subsection{Bypasses in overtwisted contact structures} 
\begin{figure}[htb] 
\begin{center}
\includegraphics{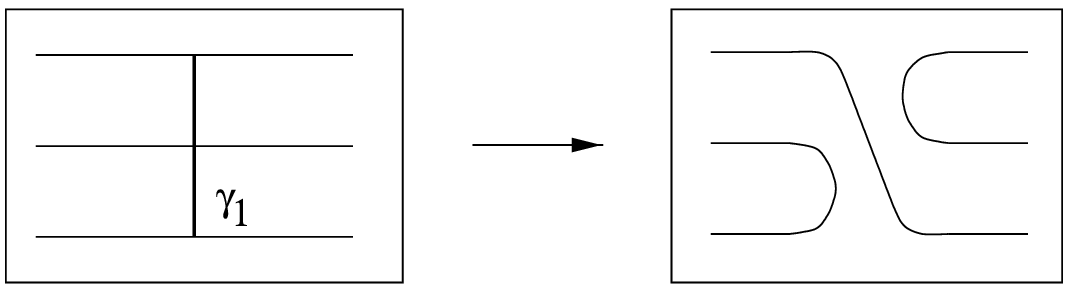}
\end{center}
\caption{\label{b:bypass}}
\end{figure}

\mlabel{s:bypasses}
Consider a contact manifold $(N,\CC)$ and a convex surface $\Sigma\subset N$ which has empty or Legendrian boundary. Let $\Gamma_{\Sigma}$ be a dividing set of $\Sigma$.
\begin{defn} \mlabel{d:bypasses}
A {\em bypass} for $\Sigma$ is an embedded half--disc $D$ in $N$ whose singular foliation has the following properties:
\begin{itemize}
\item[(i)] $\partial D$ is the union of two Legendrian arcs $\gamma_1,\gamma_2$ which intersect at their endpoints.
\item[(ii)] $D$ intersects $\Sigma$ transversally along $\gamma_1$. There are no other intersection points.
\item[(iii)] $D$ admits an orientation such that  
\begin{itemize}
\item there are exactly two positive singularities along $\gamma_1$. These are the endpoints of $\gamma_1$. They are elliptic.
\item there is exactly one negative singularity on $\gamma_1$. It is elliptic.
\item there are only positive singularities along $\gamma_2$. They alternate between elliptic and hyperbolic.
\end{itemize}
\item[(iv)] $\gamma_1$ intersects $\Gamma_\Sigma$ in exactly three points. The intersections are transverse and correspond to the singularities of the singular foliation on $D$ along $\gamma_1$.
\item[(v)] The dividing set of $D$ has exactly one connected component.
\end{itemize}
\end{defn}
A bypass for $\Sigma$ allows us to isotope $\Sigma$ such that the isotopy type of the dividing set changes as in \lemref{l:attach}. They where introduced by K.~Honda in \cite{honda}. Requirement (v) in \defref{d:bypasses} does not appear in \cite{honda} since all contact structures considered in that article are tight. In this situation, the dividing set $\Gamma_D$ of $D$ is determined (up to isotopy) by (i)--(iv) while this is not true for overtwisted contact structures. Requirement (v) is necessary for the following bypass attachment lemma. 
\begin{lem}[Honda, \cite{honda}] \mlabel{l:attach}
Let $D$ be a bypass for a convex surface $\Sigma$. There exists a neighborhood of $\Sigma\cup D$ which is diffeomorphic to $\Sigma\x[0,1]$ such that $\Sigma\x\{i\}$ is convex for $i=0,1$. The dividing set of $\Sigma\x\{1\}$ can be obtained from the dividing set of $\Sigma\x\{0\}$ as in \figref{b:bypass}. (The bypass is attached to the front. \figref{b:bypass} represents only a neighborhood of the attaching region of $D$.)
\end{lem} 
Gluing together two bypasses along the boundary component containing the singularities with alternating signs one obtains a disc with a singular foliation like the singular foliation on an overtwisted disc. Thus one can think of a bypass as one half of an overtwisted disc and it is not surprising that it is much easier to find bypasses if the contact structure is overtwisted rather than tight.

In the following proposition we assume that $\gamma_1$ is a Legendrian arc contained in the convex surface $\Sigma$. Since every arc $\gamma_1'$ which is transverse to $\Gamma$ and intersects $\Gamma$ as $\gamma_1$ in \propref{p:bypass existence} is non--isolating, we can apply \thmref{t:LeRP} to find a $C^0$ small admissible isotopy of $\Sigma$ such that $\gamma_1'$ becomes Legendrian, i.e. the assumption that $\gamma_1$ is Legendrian is actually not restrictive.
\begin{prop} \mlabel{p:bypass existence}
Let $\Sigma$ be a convex surface in a contact manifold and $D_{ot}$ an overtwisted disc disjoint from $\Sigma$. Assume that $\gamma_1\subset\Sigma$ is an embedded Legendrian arc such that the endpoints lie on $\Gamma$ and $\gamma_1\cap\Gamma$ contains three points. Then there is a bypass for $\Sigma$ which intersects $\Sigma$ along the $\gamma_1$.
\end{prop}
\begin{proof}
Note that $\gamma_1$ is automatically transverse to $\Gamma$. Let $V$ be a contact vector field transverse to $\Sigma$. Consider the image $R$ of $\gamma_1$ under the flow $\varphi_t$ of $V$ for $0\le t\le\eps$ such that $\gamma_1=R\cap\Sigma$. The curves $\varphi_t(\gamma_1), 0\le t\le\eps$ are Legendrian since $\varphi_t$ is the flow of a contact vector field. If $p\in\gamma_1\cap\Gamma$, then the segment $\varphi_t(p), 0\le t\le\eps$ is Legendrian.

If $q_1,q_2$ are points in different components of $\gamma_1\setminus(\gamma_1\cap\Gamma)$, then the union $\Gamma_R$ of the curves $\varphi_t(q_i)$ with $0\le t\le\eps, i=1,2$ divide the singular foliation on $R$. Since $R$ has Legendrian boundary, \thmref{t:divide impera} implies that $R$ is convex.  We orient $R$ such that the singularities at the endpoints of $\gamma_1$ are positive. By \eqref{e:kanda rot}
\begin{align*}
\tb(\partial R) &=-2 & \rot(\partial R) &=1~.
\end{align*}
The idea is to perform a Legendrian connected sum of the knots $\partial R$ and $\partial D_{ot}$ together with a boundary connected sum of $R$ and $D_{ot}$. Doing this carefully enough, one obtains a bypass. Let us first explain the Legendrian connected sum of null--homologous Legendrian knots $K_1,K_2$ in a contact manifold, cf.~\cite{etnyrehonda2}.

Let $p_1\in K_1$ and $p_2\in K_2$. Choose a Legendrian curve $\gamma$ between $p_1$ and $p_2$ such that $\gamma$ does not meet the knots anywhere else and $\gamma$ is not tangent to $K_i$ in $p_i$ for $i=1,2$. By \corref{c:tube leg curve} there are coordinates $x,z,t$ on a tubular neighborhood of $\gamma$ such that $\gamma$ corresponds to the $t$--axis and the contact structure is defined by $dz-x\,dt$. After a $C^0$--small Legendrian isotopy of $K_1$ and $K_2$ we may assume that $p_1,p_2$ are cusp points of the front projection, i.e the projection to the $z,t$--plane, and that the knots are oriented as in \figref{b:confront}. The Legendrian connected sum $K_1\# K_2$ is then formed using the dashed curves in \figref{b:confront}.
\begin{figure}[ht] 
\begin{center}
\includegraphics{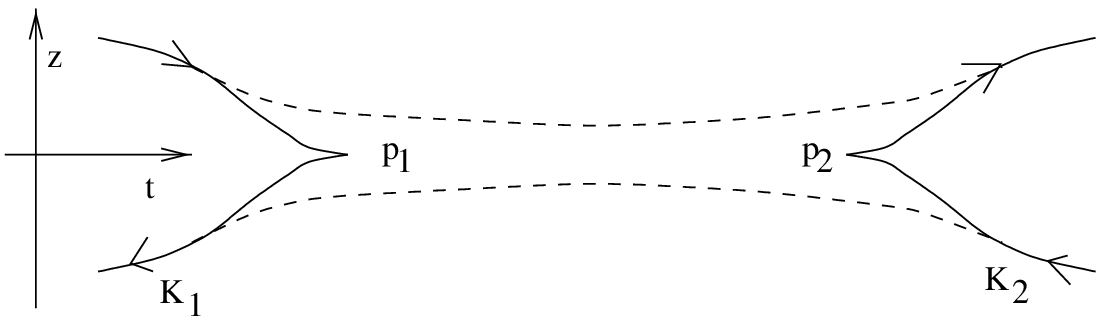}
\end{center}
\caption{\label{b:confront}}
\end{figure}
 
Now consider Seifert surfaces $\Sigma_1$ of $K_1$ and $\Sigma_2$ of $K_2$. We assume that $\Sigma_1$, respectively $\Sigma_2$, coincides with translates of $K_1$ in the negative $t$--direction, respectively of $K_2$ in the positive $t$--direction, on a neighborhood of $p_1$, respectively $p_2$. Such neighborhoods are shown in \figref{b:confront}. If we orient $\Sigma_1$ and $\Sigma_2$ such that $K_1$ and $K_2$ are oriented as boundaries, then $p_1$ is a negative singularity and $p_2$ is a positive singularity. 

We use the ribbon which is bounded by the dashed curves in \figref{b:confront} to construct a Seifert surface $\Sigma_1\#\Sigma_2$ for the knot $K_1\# K_2$. According to \cite{etnyrehonda2}
\begin{align}
\label{e:consumtb} \tb(K_1\# K_2,[\Sigma_1\#\Sigma_2]) & = \tb(K_1,[\Sigma_1]) + \tb(K_2,[\Sigma_2]) + 1 \\
\label{e:consumrot} \rot(K_1\# K_2,[\Sigma_1\#\Sigma_2]) & = \rot(K_1,[\Sigma_1]) + \rot(K_2,[\Sigma_2])~.
\end{align}
Now we apply this to the Seifert surfaces $R=\Sigma_1$ of $\partial R=K_1$ and $D_{ot}=\Sigma_2$ of $\partial D_{ot}=K_2$ where $D_{ot}$ is oriented such that $\tb(\partial D_{ot})=0$ and $\rot(\partial D_{ot}) = -1$. Then we perturb $R\# D_{ot}$ to a convex surface with boundary $\partial R\#\partial D_{ot}$.  

The difficulty in showing that the perturbed surface is a bypass is to establish (v) of \defref{d:bypasses}. In order to do this we reduce the region where the Legendrian connected sum is performed such that the contact structure on this region is tight. 

By \thmref{t:LeRP} and \thmref{t:giroux flex} we can assume that the singular foliation on $D_{ot}$ is of the form indicated in \figref{b:otdisc+} where the thickened circle is the dividing set. Then we can decompose $D_{ot}$ into two half--discs bounded by Legendrian unknots with Thurston--Bennequin invariant $-1$ and rotation number $0$. The half--discs are separated by straight Legendrian arcs. The singular foliation near the unknots is in the standard form used in \propref{p:convexity with boundary}.
\begin{figure}[htb] 
\begin{center}
\includegraphics{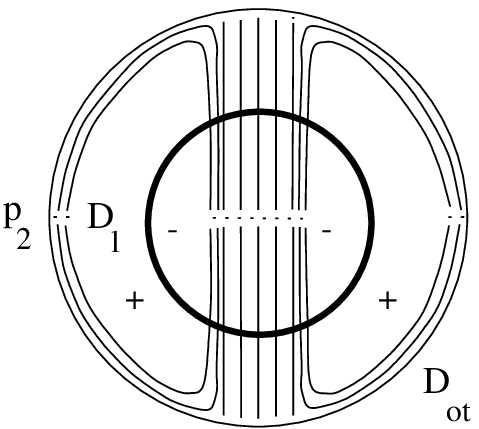}
\end{center}
\caption{\label{b:otdisc+}}
\end{figure}

By the last part of \propref{p:convexity with boundary} we can pretend that we form a boundary connected sum of the surfaces $R$ with the left part $D_l$ of $D_{ot}$ and that the perturbation of $R\# D_{ot}$ only changes the half--disc $D_l$. The presence of the Legendrian curves in the middle of $D_{ot}$ prevents an interaction between the left and the right part of $D_{ot}$. 

The union of tubular neighborhoods of $R$, of the Legendrian arc connecting $R$ with $D_{ot}$ and of the left part $D_l$ of $D_{ot}$ can be recovered in tight contact manifolds: $D_l$ can be obtained applying \thmref{t:giroux flex} to a bypass in a tight contact manifold, hence a neighborhood of $R\# D_{l}$ is tight. 

The Thurston--Bennequin invariant of $\partial R\# \partial D_l$ is $-2$ by \eqref{e:consumtb}. From \eqref{e:kanda rot} and \thmref{t:giroux tight} it follows that the dividing set on $R\# D_l$ (after this surface is perturbed to a convex surface) consists of exactly two arcs with endpoints on $\partial R\#\partial D_l$ and no closed components. The notations $R\# D_l$ and $R\# D_{ot}$ are slightly misleading because $D_l$, respectively $D_{ot}$, is {\em not} a subset of $R\# D_l$, respectively $R\# D_{ot}$, after these surfaces are smoothed and made convex.

We are left with the two possibilities for the isotopy type of the dividing set of $R\# D_{ot}$  shown in \figref{b:divbypass}. 
\begin{figure}[htb] 
\begin{center}
\includegraphics{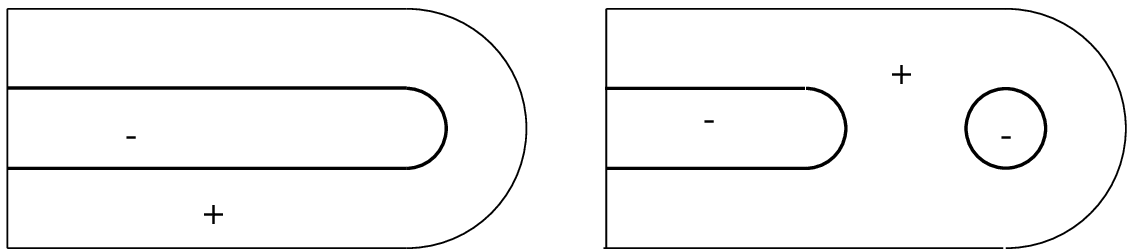}
\end{center}
\caption{\label{b:divbypass}}
\end{figure}
According to \eqref{e:kanda rot} the boundary of the left part of \figref{b:divbypass} has rotation number $1$ while in the right part its rotation number is $-2$. By \eqref{e:consumrot} the rotation number of $\partial R\# \partial D_{ot}$ is $1$. This shows that the dividing set on the Seifert surface of $\partial R\# \partial D_{ot}$ contains no closed components.

The remaining conditions (i), (ii) and (iv) in \defref{d:bypasses} are satisfied by construction. Condition (iii) can be achieved using \thmref{t:giroux flex}. Thus $R\# D_{ot}$ yields a bypass.
\end{proof}

%%%%%%%%%%%%%%%%%%%%%%%%%%%%%%%%%%%%%%%%%%%%%%%%%%%%%%%%%%%%%%%%%%%%%%%%%%
\section{Round handle decompositions and model Engel structures} \mlabel{s:round model}

This section provides the building blocks for our construction of Engel manifolds. We first recall several facts from \cite{asimov1} about round handle decompositions. In \secref{s:perturbed prolong} we describe a perturbed version of the prolongation construction and first examples of Engel structures which are compatible with a round handle decomposition of the underlying manifold. Perturbed prolongation will be used for the construction of model Engel structures on round handles in Sections \ref{s:index zero models} to \ref{ss:model 3-handles}. 

\subsection{Round handle decompositions}

Round handles were used by D.~Asimov (\cite{asimov1}) for the study of flow manifolds. A flow manifold is a manifold $M$ with a non--singular vector field $W$ transverse to the boundary which points outwards along $\partial_+M$ and inwards along $\partial_-M$. 
\begin{defn} \mlabel{d:Rk}
$R_l = D^l\x D^{n-l-1}\x S^1$ is a {\em round handle} of dimension $n$ and index $l\in\{0,\ldots,n-1\}$. The boundary $\partial R_l$ is the union of the two subsets
\begin{align*}
\partial_-R_l & = \partial D^l\x\phantom{\partial}D^{n-l-1}\x S^1 \\
\partial_+R_l & = \phantom{\partial}D^l\x\partial D^{n-l-1}\x S^1~.
\end{align*}
\end{defn}
A round handle is attached to a flow manifold using an embedding $\psi : \partial_-R_l \lra \partial_+M$. One can extend $W$ to a non--singular vector field on $M\cup R_l$ which is again transverse to the boundary. This is impossible if one attaches an ordinary handle $h_l=D^l\x D^{n-l}$ of index $l\in\{0,\ldots,n\}$ to $M$ since this changes the Euler characteristic by $(-1)^l$. 

Now assume that $W$ spans the characteristic foliation of an Engel structure on $M$. If one extends the Engel structure from $M$ to $M\cup R_l$ it is clear from the properties of transverse hypersurfaces that it is useful to ensure that the boundary of the new Engel manifold is again transverse to the characteristic foliation. Therefore round handles are suitable building blocks for the construction of Engel structures.
%------------------ def of a decomposition into round handles
\begin{defn} \mlabel{d:round decomp}
If $M$ is obtained from the disjoint union of finitely many round handles of index $0$ by attaching round handles of higher index successively, i.e.
\begin{equation*} 
M=\left(\ldots\left(\bigcup R_0\right)\cup_{\psi_1}R_{\beta_1}\ldots\right)\cup_{\psi_r}R_{\beta_r}
\end{equation*}
with $\beta_i\in\{1,\ldots,n-1\}$ for $i\in\{1,\ldots,r\}$, then we have a {\em round handle decomposition} of $M$. 
\end{defn}
We can rearrange a given round handle decomposition of a manifold such that the round handles are ordered according to their index. Contrary to the case of ordinary handles, two round handles of the same index can not be interchanged in general.
 
If a closed manifold $M$ admits a round handle decomposition, then its Euler characteristic $\chi(M)$ has to vanish because we can use the round handle decomposition to find a non--singular vector field on $M$. The converse direction is covered by the following theorem. 
\begin{thm}[Asimov, \cite{asimov1}] \mlabel{t:round decomp euler 0}
A closed, connected manifold of dimension $n\neq 3$ admits a decomposition into round handles if and only if $\chi(M)=0$. This decomposition can be chosen such that there is only one round $0$--handle and one round $(n-1)$--handle. 
\end{thm}
The analogous statement is wrong for manifolds of dimension $3$, cf.~\cite{morgan}. In order to prove the last part of \thmref{t:round decomp euler 0} one follows the proof in \cite{asimov1} starting with a decomposition of $M$ into ordinary handles with only one handle of index $0$, respectively $n$.

One can decompose $R_l$ into one ordinary handle of index $l$ and another one of index $l+1$ as follows
\begin{align} \label{e:decomp Rk}
\begin{split}
R_l & = D^l\x D^{n-l-1}\x S^1 = D^l\x D^{n-l-1}\x (D^1\cup D^1) \\
    & = \left(D^l\x (D^{n-l-1}\x D^1)\right)\cup \left((D^l\x D^1)\x D^{n-l-1}\right) \\
    & = h_l \cup h_{l+1}~.
\end{split}
\end{align}
This allows us to obtain decompositions into ordinary handles from round handle decompositions.

\subsection{Perturbed prolongation} \mlabel{s:perturbed prolong}
From now on we consider only $4$--dimensional round handles. For the construction of Engel structures we will fix a set of particular Engel structures on $R_l$ for each $l=0,1,2,3$.
\begin{defn}
A model Engel structure on a round handle of index $l$ is an Engel structure on $R_l$ such that the characteristic foliation is oriented and  transverse to $\partial_+R_l$ and $\partial_-R_l$. It points outwards along $\partial_+R_l$ and inwards along $\partial_-R_l$. 
\end{defn}
The usual prolongation construction is described in \exref{e:prolong} and deformations of certain prolonged Engel structures are discussed in \cite{montgomery}. We now describe perturbations of Engel structures which we will use later for the explicit construction of model Engel structures.
 
Let $\alpha$ be a contact form on a compact manifold $N$ such that there is a contact vector field $V$ for $\CC=\ker(\alpha)$ which is transverse to $\partial N$. Moreover we assume that $\CC$ is trivial as a bundle and we fix a trivialization $C_1,C_2$ of $\CC$. We use the same notation for the horizontal lifts to $N\x S^1$. 
\begin{prop} \mlabel{p:perturbed prolong}
For $k\in\Z$ the distribution $\DD_k$ on $N\x S^1$ spanned by 
\begin{equation*} %\label{e:perturbed prolong}
W=\partial_{t} + \eps V \textrm{ and } X_k=\cos(kt)C_1+ \sin(kt)C_2
\end{equation*}
is an Engel structure on $N\x S^1$ for $\eps>0$ small enough and $k\neq 0$. For $k=0$ we obtain an Engel structure if and only if $C_1$ and $[V,C_1]$ are linearly independent everywhere. 

Then the characteristic foliation of $\DD_k$ is spanned by $W$. In particular, the boundary is transverse and both $\CC$ and the contact structure on $\partial N\x S^1$ induce the same singular foliation on $\partial N\x \{1\}$.
\end{prop}
\begin{proof}
In order to show that $\DD_k$ is an Engel structure we first calculate
\begin{align} \label{e:perturbed prolong 1}
\begin{split}
[W,X_k] & = -k\sin(kt)C_1 + k\cos(kt)C_2 \\
        & \quad + \eps(\cos(kt)[V,C_1]+\sin(kt)[V,C_2])~.
\end{split}
\end{align}
If $k\neq 0$ and $\eps>0$ is small enough, then $W,X_k$ and $[W,X_k]$ are linearly independent. The case $k=0$ is obvious. Since $V$ is a contact vector field $[\DD_k,\DD_k]$ is generated by $W,C_1$ and $C_2$. Hence $\EE=[\DD_k,\DD_k]$ is defined by 
\begin{equation} \label{e:perturbed prolong EE}
\beta=\alpha-\eps\alpha(V)\,dt~.
\end{equation}
From the fact that $\alpha$ is a contact form it follows that $\beta$ defines an even contact structure for all $\eps>0$. This shows that $\DD_k$ is an Engel structure. 

Notice that $\beta(\eps V+\partial_t)=0$. Let $h$ be the function with the property $L_V\alpha=h\alpha$. Then
\begin{align*}
L_{(\eps V+\partial_t)}(\alpha-\eps\alpha(V) dt) & = \eps L_V\alpha - \eps^2(L_V\alpha(V))dt \\
 & = \eps h(\alpha-\eps\alpha(V))
\end{align*} 
implies that $W$ spans the characteristic foliation of $\DD_k$.
\end{proof}
If $N=h_l$, then we obtain a model Engel structure on $R_l=h_l\x S^1$ if $V$ points outwards along $\partial_+h_l=D^l\x \partial D^{3-l}$ and inwards along $\partial_-h_l=\partial D^l\x D^{3-l}$.
  
Let us describe the orientation of $\EE=[\DD_k,\DD_k]$ : If $k=0$ we obtain an Engel structure only if $C_1$ and $[V,C_1]$ are linearly independent and then $\EE$ is oriented by $W,C_1,[V,C_1]$. For $k>0$ the even contact structure $\EE$ is oriented by $W,C_1,C_2$. If $k<0$, then $\EE$ has the opposite orientation. 

If $p\in\partial N$ lies on the characteristic surface $\Sigma=\{\alpha(V)=0\}$ of $V$, then $\{p\}\x S^1$ is tangent to the contact structure on $\partial N\x S^1$. For later use we determine the rotation number along these Legendrian curves. If $k=0$ the rotation number is zero since the intersection line field along $\partial N\x S^1$ is $S^1$--invariant. For $k\neq 0$ we may assume that $C_1=V$ at $p$. Clearly
\begin{equation} \label{e:before projection}
V=C_1 = \cos(kt)X_k - 1/k\sin(kt)[W,X_k] + R_\eps
\end{equation}
where the correction term $R_\eps$ satisfies $\lim_{\eps\to 0}R_\eps=0$. If we project all vector fields in \eqref{e:before projection} along $W$ to $\partial N\x S^1$ we obtain an analogous expression for $\partial_t$. Notice that the projection of $X_k$ to $\partial N\x S^1$ spans the intersection line field. Hence the rotation number along $\{p\}\x S^1$ is $-|k|$ (cf.~\defref{d:rotzahl}). 

In order to obtain model Engel structures with positive rotation number along $\{p\}\x S^1$ one can replace $W$ by $-\partial_t+\eps V$, but doing so one changes the contact structure on $\partial N\x S^1$. In the cases $N=h_l$ for $l=0,\ldots,3$ considered later we will construct a diffeomorphism $f$ of $R_l=h_l\x S^1$ which preserves the contact structure on large parts of $\partial_-R_l$ and reverses the orientation of $\{p\}\x S^1$. Pushing forward the model Engel structures with $f$ one obtains model Engel structures with non--negative rotation numbers along $\{p\}\x S^1$. 

Before we discuss explicit model Engel structures let us first describe examples of Engel manifolds with a decomposition into round handles which carry model Engel structures.
 
If $N$ is decomposed into ordinary handles we obtain a round handle decomposition of $N\x S^1$ such that the round handles are products of $S^1$ with ordinary handles from the decomposition of $N$. If $V$ is transverse to the boundaries of all handles and points into the right directions, then $N\x S^1$ carries an Engel structure such that all round handles carry model Engel structures. This situation arises in the context of convex contact structures defined in \cite{eliashberggromov}.
\begin{defn}  \mlabel{d:convcontstr}
A contact structure $\CC$ on a manifold $N$ is {\em convex} if there is a proper Morse function $g : N\rightarrow [0,\infty)$ and a complete contact vector field $V$ which is a pseudo--gradient for $g$, i.e. there is a Riemannian metric and a positive function $s$ such that $L_Vg\ge s\|dg\|^2$.
\end{defn}
According to \cite{giroux icm} every contact structure on a $3$--manifold is convex. From $g$ one obtains a handle decomposition of $N$ and applying \propref{p:perturbed prolong} to convex contact structures one obtains examples of Engel structures with adapted round handle decomposition.

\subsection{Model Engel structures on round handles of index $0$} \mlabel{s:index zero models}

For our proof of \thmref{T:EXIST} it is useful to have overtwisted contact structures on transverse boundaries of Engel manifolds. Since we will start our construction of Engel manifolds with a round handle of index $0$ we define the model Engel structures on $R_0$ such that the contact structure on $\partial_+R_0$ is overtwisted. 

We use cylindrical coordinates on $\R^3$. Consider the overtwisted contact structure $\ker(\cos(r^2)\,dz+\sin(r^2)\,d\varphi)$ and the ball of radius $r_0=3\pi/2$. The singular foliation on the boundary $S(r_0)$ admits a dividing set with three connected components. From \thmref{t:divide impera} it follows that $S(r_0)$ is convex, i.e. there is a contact vector field $V$ transverse to $S(r_0)$ which points out of the ball. Since the ball is contractible the contact structure admits a trivialization $C_1,C_2$. 

By \propref{p:perturbed prolong} we obtain model Engel structures $\DD_k, k\in\Z\setminus\{0\}$ on $R_0$. At this point we do not need to deal with the case $k=0$. The singular foliation on $S(r_0)\subset\partial_+R_0$ is the same as the original singular foliation on $S(r_0)$. By \thmref{t:giroux tight} the contact structure on $\partial_+R_0$ is overtwisted. 

Since $\pi_1(\SO(4))=\Z_2$ there are exactly two homotopy classes of oriented framings on $R_0$. The homotopy class of the Engel framing of $\DD_k$ depends only on the parity of $k$. So we have shown
\begin{prop} \mlabel{p:index 0 models}
For each of the two homotopy classes of oriented framings on $R_0$ there is a model Engel structure on $R_0$ such that the Engel framing is homotopic to the given framing of the tangent bundle of $R_0$. Moreover, the induced contact structure on the boundary is overtwisted.
\end{prop}  
\begin{rem}
J.~Adachi asks in \cite{adachi} whether it is possible to find an Engel structure on $N\x[0,1]$ such that the boundary is transverse and the induced contact structures $\CC_0$ on $N\x\{0\}$ and $\CC_1$ on $N\x\{1\}$ are equivalent to given contact structures on $N$. One can show that $\CC_0$ and $\CC_1$ have to be homotopic as plane fields on $N$ as in \lemref{l:unique homotopy class of plane fields}. The following example shows that $\CC_0$ and $\CC_1$ do not have to be equivalent.

Consider the contact structure from above on $B(r_0)$. We remove a ball contained in a Darboux chart with convex boundary $S'$ from the interior of $B(r_0)$ and choose the contact vector field $V$ such that it is transverse to both $S'$ and $S(r_0)$. By \propref{p:perturbed prolong} we obtain an Engel structure on a manifold diffeomorphic to $S^2\x [0,1]\x S^1$ such that 
\begin{itemize}
\item the characteristic foliation is transverse to the boundary
\item the contact structure on one boundary component is tight while it is overtwisted on the other boundary component.
\end{itemize}
This answers the question completely for $N=S^2\x S^1$ by the classification of contact structures on $S^2\x S^1$. 
\end{rem}

\subsection{Model Engel structures on round handles of index $1$} \mlabel{s:index one models}

On $R_1=D^1\x D^2\x S^1$ we denote the coordinate on $D^1$ by $x$, the coordinates on $D^2$ by $y_1, y_2$ and the coordinate on $S^1$ by $t$. On $h_1=D^1\x D^2$ consider 
\begin{align}
\begin{split} \label{e:runde 1-henkel}
\CC&=\ker(\alpha=-dy_1 -y_2\,dx -1/2\,x\,dy_2) \\ 
V&=1/2\,y_1\rve{y_1} + y_2\rve{y_2}-1/2\,x\rve{x}
\end{split}
\end{align}
These choices are motivated by \cite{weinstein, TopcharStein}: The resulting even contact structure on $R_1$ is defined by the contraction of the symplectic form $dy_1\ww dt+dx\ww dy_2$ with the Liouville vector field $\partial_t+V$. Notice that $V$ is transverse to $\partial_\pm h_1$ and points outwards along $\partial_+h_1=D^1\x \partial D^2$ and inwards along $\partial_-h_1=\partial D^1\x D^2$. We apply \propref{p:perturbed prolong} to
\begin{align} \label{e:runde 1-henkel 2}
\begin{split}
C_1&=\rve{y_2}-\rve{x}+(y_2-1/2\,x)\rve{y_1} \\
C_2&=[V,C_1]=-\rve{y_2}-1/2\rve{x}+1/2(y_2+x)\rve{y_1}~.
\end{split}
\end{align}
with $\eps=1$. The properties of the resulting distributions $\DD_k$ are summarized in the following proposition. 
\begin{prop} \mlabel{p:1 modelle}
The distributions $\DD_k$ are model Engel structures on $R_1$ for all $k\in\Z$. They have the following properties. 
\begin{itemize}
\item[(i)] The even contact structure $\EE=[\DD_k,\DD_k]$ is defined by $\beta=\alpha+1/2\,y_1\,dt$. The orientation of the contact structure on $\partial_-R_1$ and $\partial_+R_1$ with respect to the restriction of $d\beta$ is positive if $k\ge 0$ and negative if $k<0$ .
\item[(ii)] 
The curves $\gamma_\pm=\{\pm 1\}\x\{0\}\x S^1$ are Legendrian. The rotation number along them is $-|k|$.
\item[(iii)]
The rotation number of the intersection line field with respect to the Legendrian vector field
$$
Z=y_2\rve{t} + 1/2\,y_1\rve{x}
$$
along $\{0\}\x\{y_1=0, y_2=1\}\x S^1$ is $-|k|$ and it equals $0$ along $\{0\}\x S^1\x\{0\}$.
\end{itemize}
\end{prop}
This follows from the discussion of \propref{p:perturbed prolong} and calculations using \eqref{e:runde 1-henkel} and \eqref{e:runde 1-henkel 2}. 
\begin{rem} \mlabel{r:contact contraction} 
The contact structure on $\partial_-R_1$ is defined by $\beta_-=-dy_1+1/2\,y_1\,dt-1/2\,x\,dy_2$ with $x=\pm 1$. It is invariant under the map $f_s$ induced by $(y_1,y_2)\lmt (sy_1,sy_2)$ for $s\in(0,1)$. In particular, if $\psi : \partial_-R_1\lra \partial_+M$ is an embedding preserving contact structures, then the same is true for $\psi\circ f_s$. We can choose $s$ so small that the image of $\psi\circ f_s$ is contained in a given tubular neighborhood of $\psi(\gamma_\pm)$.  
\end{rem}

Up to now we have a model Engel structure on $R_1$ for each orientation of the contact structure on $\partial_-R_1$ and each homotopy class of intersection line fields with negative rotation number along $\gamma_\pm$. Moreover, we have a model Engel structure with rotation number $0$ along $\gamma_\pm$ for one of the two possible orientations of the contact structure on $\partial_-R_1$.

In order to obtain the missing possibilities we extend the contact diffeomorphism 
$$
f : (x,y_1,y_2,t) \lmt (x,y_1,-y_2-4xy_1,-t)
$$
of a neighborhood of $\gamma_\pm$ in $\partial_-R_1$ to a diffeomorphism of $R_1$. This map preserves $\gamma_\pm$ and the contact structure near $\partial_-R_1$ but it reverses the orientations of $\gamma_\pm$ and of the contact structure. We push forward the Engel structures $\DD_k$ we have obtained so far using the extension of $f$. Then the rotation number along $\gamma_\pm$ changes its sign. 

Up to now all model Engel structures induce the same contact framing along $\gamma_\pm$. In order to realize other contact framings one can push forward the model Engel structures obtained so far with the self diffeomorphisms $\Theta_m$ of $R_1$
\begin{align} \mlabel{e:Theta}
\begin{split}
\Theta_m :  D^2\x I\x S^1 & \lra D^2\x I\x S^1 \\
           (y_1,y_2,x,t)  & \lmt \begin{array}{l} (\cos(mt)y_1-\sin(mt)y_2, \\ \quad \sin(mt)y_1+\cos(mt)y_2,x,t) \end{array}
\end{split}
\end{align}
with $m\in\Z$. This map preserves $\gamma_\pm$ as well as $\partial_-R_1$. Now we have shown the following proposition.
\begin{prop}
For each oriented framing of $\gamma_\pm\subset\partial_-R_1$ and $k\in\Z$ there is a model Engel structure on $R_1$ such that the contact framing is homotopic to the given framing and the rotation number along $\gamma_\pm$ is $k$.

Given one such model Engel structure, there is another model Engel structure which induces the same contact structure on $\partial_-R_1$ with the same rotation number along $\gamma_\pm$ but which induces the opposite orientation of the contact structure on $\partial_-R_1$.   
\end{prop}
For $k,m\in\Z$ let $\DD_{k,m}$ for $k,m\in\Z$ be the push forward with $\Theta_m$ of a model Engel structure with rotation number $k$ along $\gamma_\pm$. The orientation of the contact structure on $\partial_\pm R_1$ does not appear in this notation.

\subsection{Model Engel structures on round handles of index $2$} \mlabel{s:modelle index 2}

The model Engel structures on $R_2=D^2\x D^1\x S^1$ depend on two parameters $n,k$ which correspond to the homotopy class of the intersection line field near the torus $T^2_0=\partial D^2\x\{0\}\x S^1\subset\partial_-R_2$. The contact structure on $\partial_-R_2$ is essentially independent of the model Engel structure.   
\begin{prop} \mlabel{p:zwei henkel modell}
Given integers $n\in\Z$ and $k\in\Z\setminus\{0\}$ there is a model Engel structure $\DD=\DD_{k,n}$ on $R_2$ with the following properties.
\begin{itemize}
\item[(i)] The characteristic foliation of $\DD$ can be oriented such that it points outwards along $\partial_+R_2$ and inwards along $\partial_-R_2$.
\item[(ii)] The singular foliation on $T^2_0=\partial D^2\x\{0\}\x S^1\subset\partial_-R_2$ is in standard form. In particular, $T^2_0$ is convex. The Legendrian ruling corresponds to the first factor of $T_0^2=\partial D^2\x\{0\}\x S^1$. The dividing curves are tangent to the last factor. 
\item[(iii)] The rotation number of the intersection line field along $\gamma=\partial D^2\x\{0\}\x\{1\}$ (with its orientation as boundary $\partial D^2$) is $2n$. 
\item[(iv)] The rotation number of the intersection line field along the Legendrian divides (with the canonical orientation of the last factor of $\partial D^2\x\{0\}\x S^1$) is $k\neq 0$. 
\item[(v)] The orientation of the contact structure on $\partial_+R_2$ can be chosen freely. 
\end{itemize}
All model Engel structures induce the same contact structure on a neighborhood of $T^2_0\subset\partial_-R_2$.
\end{prop}
\begin{proof} 
We choose the even contact structure on $R_2=D^2\x D^1\x S^1$ first. The starting point is a singular foliation $\FF$ on a disc $D^2$ where we use polar coordinates $(r,\varphi)$. The coordinate on $D^1$ is $x$. 

On $A=\{r>1/2\}$ we define $\FF$ by $\cos(\varphi)dr$. In addition we require that $\FF$ admits a dividing set $\Gamma$ with the following properties.
\begin{itemize}
\item[(i)] The straight arc $\gamma_0$ from $(r=1,\varphi=0)$ to $(r=1,\varphi=\pi)$ lies in $\Gamma$. In addition, $\Gamma$ contains $|n|$ closed circles.
\item[(ii)] Except for $\gamma_0$, every component of $\Gamma$ is closed and bounds a disc containing no other components of $\Gamma$. If $n>0$ all closed components lie below $\gamma_0$ in $D^2\setminus\gamma_0$ and if $n<0$ they lie in the upper half--disc.
\end{itemize}
\figref{b:r2} shows a possible $\FF$ for $n=2$. The thickened curves divide $\FF$. Similar singular foliations can be found other $n$. 
\begin{figure}[htb] 
\begin{center}
\includegraphics{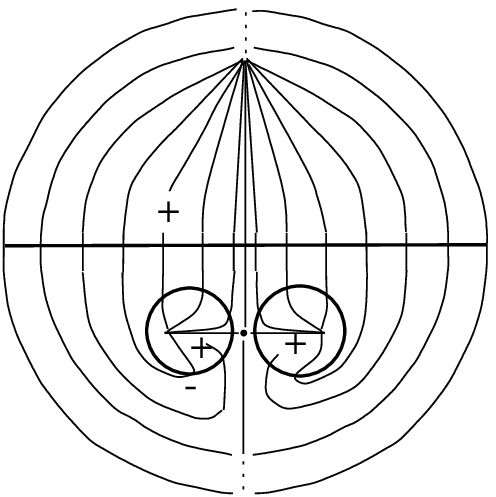}
\caption{\label{b:r2}}
\end{center}
\end{figure}
By \thmref{t:divide impera} there is an $\R$--invariant positive contact structure $\CC$ on $D^2\x\R$ such that the induced singular foliation on $D^2\x\{0\}$ is $\FF$. We choose an $\R$--invariant contact form $\alpha$ for $\CC$. such that 
$$
\alpha= \cos(\varphi)\,dr+\sin(\varphi)\,dx
$$ 
on $A\x R$. The coordinate corresponding to the $\R$--factor is $x$. This choice fixes an orientation of the contact structure. In order to find a contact vector field $V$ and a $2$--handle $h_2\subset D^2\x\R$ such that $V$ is transverse to $\partial h_2$ we need to take some care since we know nothing about the region $r<1/2$ except that $\partial_x$ is a contact vector field everywhere. We focus first on $A\x\R$. 

Let $g_1,g_2$ be smooth functions on $A\x\R$ depending only on $x$. By the proof of \propref{p:convectfunctions} 
\begin{equation*} 
V = g_1(x)\rve{r} -\left(g_1'(x)\cos^2(\varphi)+g_2'(x)\sin(\varphi)\cos(\varphi)\right)\rve{\varphi} + g_2(x)\rve{x}
\end{equation*}
is the contact vector field associated to the function $h=g_1(x)\cos(\varphi)+g_2(x)\sin(\varphi).$

We choose the functions $g_1,g_2$ such that
\begin{align}
\nonumber
g_1(x) & =\left\{ \begin{array}{ll} 0   & \textrm{ for } |x|\ge 1 \\
                                   -1   & \textrm{ for } |x|\le\frac{3}{4} \end{array} \right. & \hspace{3mm}
\label{e:g2 r2}
g_2(x) & =\left\{ \begin{array}{ll}  a   & \textrm{ for } x\ge  \frac{3}{4}\\
                                     -a  & \textrm{ for } x\le -\frac{3}{4}\\
                                     0   & \textrm{ for } -\frac{1}{2}\le x\le \frac{1}{2}~. \end{array} \right.
\end{align}
for $a>0$. For these choices the contact vector field $V$ associated to $h$ can be extended by $a\cdot\mathrm{sgn}(x)\partial_x$ on $|x|\ge 1$ to a smooth contact vector field which we still denote by $V$. Finally, we extend $V$ to a contact vector field on the whole of $D^2\x\R$. (For this it is enough $\alpha(V)$ and to apply \propref{p:convectfunctions}.) The resulting contact vector field is transverse to $\partial D^2\x [-3/4,3/4]$ and points inwards there. Now consider the pair of hypersurfaces defined by the equation
\begin{align*}
|x| =5/4 - r^2/2~.
\end{align*}
Since $r\le 1$, both are contained in the region $|x|\ge 3/4$ where $g_2=\pm a$ depending on the sign of $x$.

If we fix $a$ big enough, $V$ is transverse to the hypersurfaces $\{|x|=5/4-r^2/2\}$ and it points outwards. Then
$$
h_2=\left\{(r,\varphi,x)\big|~ |x|\le 5/4 - r^2/2 \right\}
$$
is diffeomorphic to an ordinary handle of index $2$ and $V$ is transverse to both boundary components, cf.~\figref{b:henkel}. Moreover $V$ has the desired orientations along $\partial_\pm h_2$.
\begin{figure}[htb]
\begin{center}
\includegraphics{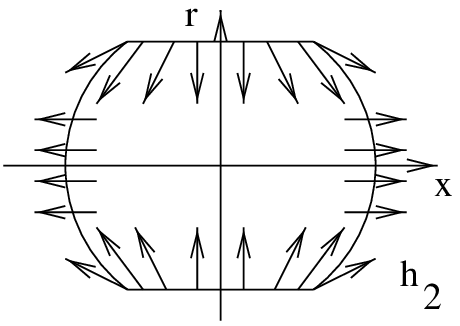}
\caption{\label{b:henkel}}
\end{center}
\end{figure}

We chose $\FF$ and $\CC$ such that $\partial D^2\subset h_2$ is Legendrian. Using \eqref{e:kanda rot} we can determine the rotation number of $\partial D^2$. We obtain 
\begin{equation} \label{e:rotzahl 2-henkel}
\rot_\CC(\partial D^2) = 2n~.
\end{equation} 
By \remref{r:rotzahlvergleich} this is (up to sign) the rotation number of any model Engel structure on $R_2$ whose even contact structure induces the singular foliation $\FF$ on $D^2\subset h^2\x S^1$. 

Now fix an oriented trivialization $C_1,C_2$. We denote the horizontal lifts of $C_1,C_2,V$ on $R_2=h_2\x S^1$ by the same symbols. As usual $t$ is the coordinate on $S^1$. By \propref{p:perturbed prolong} the distribution $\DD_k$ spanned by 
\begin{align*}
W  = \rve{t}+\eps V \hspace{4mm}\textrm{ and }\hspace{4mm} X_k = \cos(kt)C_1 + \sin(kt)C_2 
\end{align*}
is a model Engel structure if $\eps>0$ is small enough and $k\in\Z\setminus\{0\}$. The characteristic foliation of $\DD_k$ is spanned by $W$. This vector field is transverse to $\partial_\pm R_2$ and it points in the desired directions. The even contact structure $\EE=[\DD_k,\DD_k]$ is defined by $\beta=\alpha - \eps\alpha(V)\,dt$. Using the expressions for $V,\alpha,h$ and our choices of $g_1,g_2$ we obtain
\begin{align*} 
\beta & =\cos(\varphi)\,dr+\sin(\varphi)\,dx +\eps\cos(\varphi)\,dt~.
\end{align*}
on $\widetilde{A}=A\x\{-1/2\le x\le 1/2\}\x S^1$. The contact structure on $\partial_-R_2$ is defined by 
\begin{equation} \label{e:beta modelle2}
\beta\eing{\partial_-R_2} = \sin(\varphi)\,dx+\eps\cos(\varphi)\,dt-\eps g_2(x)\sin(\varphi)\,dt~.
\end{equation}
Restricting $\beta$ to $T^2_0=\partial D^2\x\{0\}\x S^1$ we see that the singular foliation on $T^2_0$ is in standard form. The curves $\varphi=\pi/2$ and $\varphi=3\pi/2$ are the Legendrian divides and the Legendrian ruling is tangent to the foliation given by the first factor in $T^2_0=\partial D^2\x\{0\}\x S^1$. 

For $k>0$, the orientation of the even contact structure is $W,C_1,C_2$. This is the orientation used in \eqref{e:rotzahl 2-henkel}, hence the rotation number of the intersection line field along $\partial D^2$ is $2n$. If $k<0$ we obtain the opposite orientation $W,C_1,-C_2$ and therefore the rotation number of the intersection line field on $\partial_-R_2$ along $\partial D^2$ has now the opposite sign, i.e. $\rot(\partial D^2)=-2n$. The rotation number along the Legendrian divide is $-|k|$ by \propref{p:perturbed prolong}.

Let us summarize the properties of the model Engel structures we have obtained up to now. Recall that $\DD_k$ depends not only on $k$ but also on the choice of $\FF$ at the beginning of the proof and that $|n|$ is the number of closed components of the dividing set of $\FF$. In the following table $|n|$ is the number of closed components of components.
\begin{center}
\begin{tabular}{|c|c|c|c|}
\hline
      & \begin{tabular}{c} Orientation  \\ of $\EE/\WW$ \end{tabular} & \begin{tabular}{c}Rotation number \\ $\partial D^2\x\{0\}\x\{1\}$\end{tabular} & \begin{tabular}{c}Rotation number\\ Legendrian divides \end{tabular}\\ 
\hline
$k>0$ & $C_1,C_2$ & $2n$ &$ -|k|$ \\
$k<0$ & $C_1,-C_2$ & $-2n$ & $-|k|$ \\
\hline
\end{tabular}
\end{center}
Notice that since we fixed $V$ and $\alpha$ on $A\x\R$, the contact structure on $\partial_-R_2$ is independent from $k,n$.

The model Engel structures with positive rotation numbers along the Legendrian divides can be obtained by applying the involution
\begin{align*}
f : R_2 & \lra R_2\\
(r,\varphi,x,t) & \lmt (r,\varphi,-x,-t)~.
\end{align*}
to the model Engel structures we have obtained so far. The contact structure on $\{-1/2\le x\le 1/2\}\subset \partial_-R_2$ is preserved by $f$ but $f$ reverses the orientation of the Legendrian divides, cf.~\eqref{e:beta modelle2}. In particular, we can compare the orientations of the contact structure and the homotopy class of the intersection line fields with the corresponding properties of $\DD_k$. The model Engel structures $f_*\DD_k$ cover the cases missing in the table above (with the exception of $k=0$).
\end{proof}
Using a singular foliation whose dividing set has more non--closed components than $\FF$ at the beginning of the proof of \propref{p:zwei henkel modell} one can construct model Engel structures on $R_2$ such that the dividing set of $T^2_0$ has $4,6,\ldots$ components.
\begin{rem} \mlabel{r:contact contraction 2}
This remark is the analogue of \remref{r:contact contraction} for round handles of index $2$. Equip $R_2$ with a model Engel structure from \propref{p:zwei henkel modell}. If $\psi : \partial_-R_2 \lra \partial_+M$ is an embedding such that the image of $T_0^2$ is convex and such that $\psi$ preserves the isotopy class of the dividing set, then one can isotope $\psi$ such that it preserves the contact structure. This follows almost immediately from \thmref{t:giroux flex} and \thmref{t:giroux unique}, but notice that the contact structure on $\partial_-R_2$ is not invariant under $\partial_x$. Since it is homotopic through contact structures to a $\partial_x$--invariant contact structure, we can use Gray's theorem (\thmref{t:gray}) to circumvent this problem.
\end{rem}

\subsection{Model Engel structures on round handles of index $3$} \mlabel{ss:model 3-handles}

We now come to the construction of model Engel structures on round handles of index $3$. In the proof of \thmref{T:EXIST} we will attach $R_3$ to an Engel manifold $M$ with transverse boundary such that the contact structure on $\partial_+M$ is overtwisted. Hence the contact structure on $\partial_-R_3$ should also be overtwisted. Contrary to round handles of index $0$ we will need {\em all} possible homotopy classes of intersection line fields on $\partial_-R_3$. In particular we also need the case $k=0$ in \propref{p:perturbed prolong}. 

If it were enough to have the unique positive tight contact structure on $\partial_-R_3\simeq S^2\x S^1$ it would be easy to describe the model Engel structures explicitly. In the overtwisted case we were not able to find explicit formulas but we still can construct the model Engel structures. We use ordinary handles of dimension $3$ to construct a contact structure together with a contact vector field $V$ and a Legendrian vector field $X$ such that $X,[V,X]$ are linearly independent.

The following construction is -- up to a small modification -- an example of Giroux's construction of convex contact structures in \cite{giroux}. 

\subsubsection{Contact structures on ordinary handles of dimension $3$}
On $h_1=D^1\x D^2$ let $x$ be the coordinate on the first factor and $y_1,y_2$ the coordinates on the second factor. Consider the contact structure $\CC_1$ defined by $\alpha=dy_1+y_2\,dx$ and the contact vector field 
\begin{equation} \label{e:V}
V_1=2y_2\rve{y_2}-x\rve{x}+y_1\rve{y_1}~.
\end{equation}
A simple calculation shows that the Legendrian vector field $X_1=\partial_{y_2} + y_2\partial_{y_1}-\partial_x$ has the property that $X_1,[V_1,X_1]$ are linearly independent everywhere. The characteristic surface (i.e. the points where $V$ is tangent to the contact structure) of $V$ is the strip 
$$
S_1=\{(x,y_1,y_2)\in h_1~|~y_1=xy_2 \}~.
$$ 
The contact vector field $V$ is transverse to $\partial_\pm h_1$ and it points outwards (respectively inwards) along $\partial_+h_1=D^1\x \partial D^2$ (respectively $\partial_-h_1=\partial D^1\x D^2$). The intersection of $S$ with $\partial_-h_1$, respectively $\partial_+h_1$, consists of two arcs.

On $h_2=D^2\x D^1\simeq h_1$ we use the same contact structure as on $h_1$ and let $X_2=X_1$. The contact vector field is reversed, i.e. $V_2=-V_1$, so it points outwards along $\partial_+h_2=\partial_-h_1$ and inwards along $\partial_-h_2=\partial_+h_1$. The characteristic surface $S_2$ of $V_2$ satisfies $S_2=S_1$. 

\begin{figure}[htb]
\begin{center}
\includegraphics{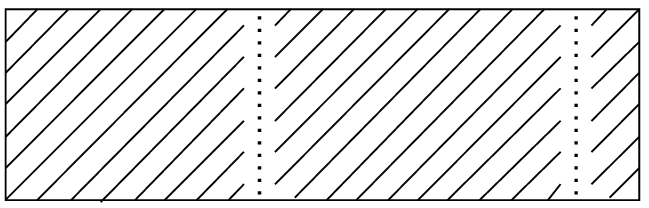}
\caption{\label{b:charaom}}
\end{center}
\end{figure}

The singular foliation on $\partial_-h_1$ consists of parallel straight segments and \figref{b:charaom} shows the singular foliation on $\partial_+h_1=\partial_-h_2$.

On $h_3=D^3$ consider the standard contact structure defined by $\alpha_3=dz+x\,dy$ and the contact vector field $V_3=-x\partial_x-2y\partial_y-3z\partial_z$. Let $X_3=\partial_x+x\partial_z-\partial_y$. Again a simple calculation shows that $[V_3,X_3]$ and $X_3$ are linearly independent everywhere. The characteristic surface of $V_3$ is 
$$
S_3= \left\{ (x,y,z)\in D_3 \big|~z=-3/2\,xy \right\}~.
$$
In particular, the dividing set $\Gamma$ on $S^2$ is connected.

Now we combine the handles with the contact vector fields and the Legendrian vector fields described above in order to construct a contact structure on $D^3$ together with a contact vector field $V$ and a Legendrian vector field $C_1$ as in \propref{p:perturbed prolong}. 

We orient $\partial_\pm h_i$ such that the orientation of $\partial_\pm h_i$ followed by $V_i$ gives the orientation of the contact structure. Moreover, we can associate a sign to each component of $\partial_\pm h_i\setminus (S_i\cap\partial_\pm h_i)$ as follows: If the orientation of $X_i,[V_i,X_i],V_i$ at $p$ is the contact orientation, then the component containing $p$ is positive. Otherwise this region is negative. Notice that along $S_i$, the vector fields $X_i,[V_i,X_i],V_i$ are not linearly independent. 
\begin{prop} \mlabel{p:D3 V X}
There is a contact structure $\CC$ on $D^3$ together with a contact vector field $V$ and a Legendrian vector field $X$ without zeroes such that $[V,X]$ and $X$ are linearly independent everywhere and $V$ is transverse to $\partial D^3$. Moreover, the dividing set on $\partial D^3$ has three connected components.
\end{prop}
\begin{proof}
Consider the overtwisted contact structure on $\R^3$ which is defined by $\alpha=\cos(r^2)\,dz+\sin(r^2)\,d\varphi$ (in polar coordinates). By \thmref{t:divide impera} the sphere $S(r_0)$ with radius $r_0=\sqrt{3\pi/2}$ around the origin is convex, i.e. there is a contact vector field $V$ transverse to $S(r_0)$. We assume that $V$ points inwards. Consider a collar $U\simeq S(r_0)\x [0,\eps], \eps>0$ of $S(r_0)$ such that $V$ corresponds to the standard vector field $\partial_\tau$ induced by the second factor. 

Clearly, one can choose a nowhere vanishing section of $\ker(\alpha)$ along $S(r_0)$ and extend it to a Legendrian vector field $X_U$ on $U$ with the desired property. Since $S(r_0)$ is simply connected, all nowhere vanishing sections of $\CC$ along $S(r_0)$ are homotopic. We orient the boundary of $U$ such that the orientation of the boundary followed by $V$ gives the orientation of the contact structure.

In order to prove the proposition we want to extend the contact structure, $X_U$ and $\partial_\tau$ to the interior of the sphere. To do so we attach the handles $h_1,h_2,h_3$. 

The dividing set $\Gamma_U$ on $\partial_+U=S(r_0)\x\{\eps\}$ consists of three connected components $\gamma_0,\gamma_+,\gamma_-$, as shown in the left part of \figref{b:surg}. Regions of $\partial_+U$ where $X,[V,X],V$ is the contact orientation are called positive and if $[X,[V,X],V$ is he opposite orientation they are called negative.
\begin{figure}[htb]
\begin{center}
\includegraphics[width=0.96\textwidth]{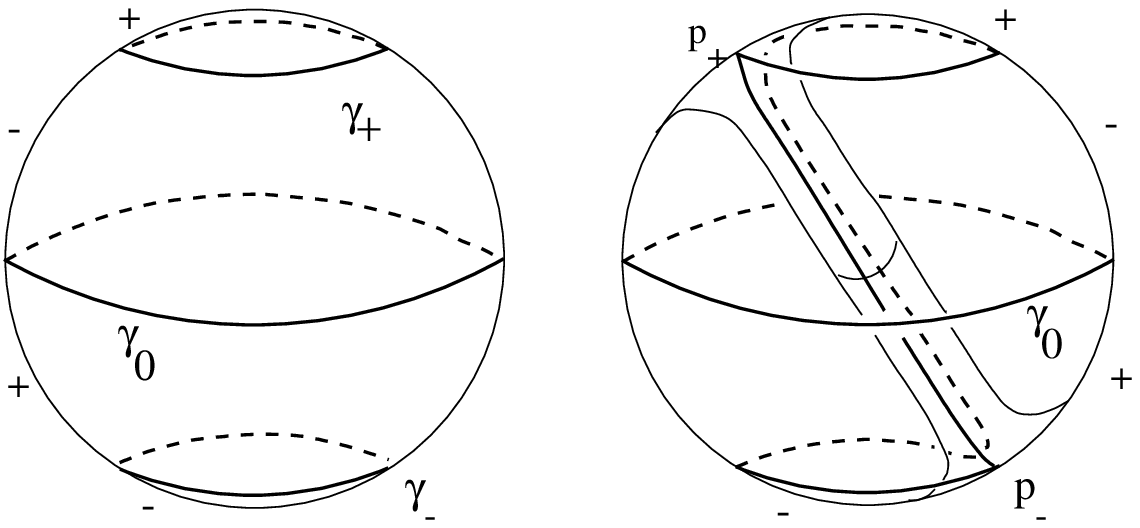}
\caption{\label{b:surg}}
\end{center}
\end{figure}

We attach $h_1$ to $U$ using an orientation preserving embedding $\psi_1 : \partial_-h_1 \lra \partial_+U$ such that $\psi_1(p_\pm)\in\gamma_\pm$ and such that $\psi_1$ preserves the singular foliation. Moreover, we assume that $\psi_1$ maps the positive, respectively negative, regions of $\partial_-h_1$ to the positive, respectively negative, regions of $\partial_+U$. Such a map exists because the dividing set on $\partial_-U$ is transverse to the singular foliation. It maps points where $V_1$ is tangent to the contact structure to $\Gamma_U$. 

Following \cite{giroux} one obtains a smooth contact structure $\CC'$ and a smooth contact vector field on $U'=U\cup_{\psi_1} h_1$ and we can modify the boundary of $U$ before we attach $h_1$ such that we obtain a smooth section $X'$ of $\CC'$ with the property that $X',[V',X']$ are linearly independent. Here we use that $\psi_1$ preserves positive and negative regions. 

The dividing set $\Gamma'$ on $\partial_+U'$ contains $\gamma_0$. The other component consists of $S\cap\partial_+h_1$ together with the parts of $\gamma_+,\gamma_-$ which do not lie in the attaching region of $h_1$. The right part of \figref{b:surg} shows $\Gamma'$. In particular $\Gamma'$ has two connected components. The homotopy class of $X'$ as a nowhere vanishing section of $\CC'$ along the torus $\partial_+U'$ is not fully determined by this construction. The restriction of $X'$ to $\gamma_0$ depends only on $X_U$ and it is therefore uniquely determined up to homotopy. On the other hand, the homotopy class of $X'$ along a curve which intersects the belt sphere of $h_1$ in exactly one point depends on the vertical modification of the boundary, cf. \remref{r:suitable vert mod 1}.

We will use this flexibility when we attach $h_2$. Using \thmref{t:divide impera} we first deform the torus $\partial_+U'$ by an admissible isotopy such that it is in standard form. The deformed space is still denoted by $U'$. We may vary the slope of the Legendrian ruling such that there is an orientation preserving embedding $\psi_2 : \partial_-h_2 \lra \partial_+U'$ which preserves singular foliations and which maps positive, respectively negative, regions of $\partial_-h_2$ to positive, respectively negative, regions of the boundary. 

Now we can extend the contact structure $\CC'$ and the contact vector field $V'$ from $U'$ to $U''=U'\cup_{\psi_2} h_2$ using $\CC_2$ and $V_2$. We denote the extensions by $\CC''$ and $V''$. When we want to extend $X'$ to $U''$ using $X_2$ we have to ensure that $\psi_{2*}(X_2)$ is homotopic to $X'$ and that the orientations of the contact structures given by $X',[V',X']$ and $X_2,[V_2,X_2]$ are coherent. 

The latter requirement is satisfied since $\psi_2$ maps positive, respectively negative, regions of $\partial_-h_2$ to positive, respectively negative, regions of $\partial_+U'$.  

In order to ensure that $\psi_2(X_2)$ is homotopic to $X'$ we use the fact that the homotopy class of $X'$ along $\partial_+U'$ is not uniquely determined. Choosing different values for $k$ in \remref{r:suitable vert mod 1} on the two connected components of the attaching region of $h_1$ we can vary the homotopy class of $X'$ along closed curves which intersect the belt sphere of $h_1$ exactly once. In this way we can ensure the $\psi_{2*}(X_2)$ is homotopic to $X'$ along the attaching curve of $h_2$. Now we can use vertical modifications of the boundary and we obtain a Legendrian vector field $X''$ on $U''$ such that $X''$ and $[V'',X'']$ are linearly independent. 

The attaching curve of $h_2$ intersects each component of $\Gamma'$ exactly once. The dividing set $\Gamma''$ of $\partial_+U''$ contains both components of $\Gamma'$ with the segments contained in the image of $\psi_2$ removed. The endpoints of the remaining curves are connected by the two segments $S\cap\partial_+h_2$. Hence $\Gamma''$ is connected and $\partial_+U''$ is diffeomorphic to a sphere. All nowhere vanishing Legendrian vector fields along $\partial_+U''$ are homotopic.

Using \thmref{t:divide impera} again, we change the singular foliation on $\partial_+U''$ such that it becomes equivalent to the singular foliation on $\partial_-h_3$. We attach $h_3$ using an embedding of $\partial_-h_3$ into $\partial_+U''$ which satisfies the same conditions on the orientations as $\psi_1$ and $\psi_2$. After a vertical modification of the boundary we can attach $h_3$ to $U''$ such that $\CC'',V'',X''$ extend to a contact structure $\CC$ on a ball, a contact vector field $V$ and a nowhere vanishing section $X$ of $\CC$ with the desired properties.
\end{proof}

\subsubsection{Model Engel structures on $R_3$}
We apply \propref{p:perturbed prolong} to $\CC,V$ and $C_1=X$ constructed in the proof of \propref{p:D3 V X}. This yields model Engel structures $\DD_k, k\in\Z$, on $R_3$ such that the contact structure on $\partial_-R_3=\partial_-U\x S^1$ is overtwisted by \thmref{t:giroux tight}. 

Because $S(r_0)$ and the contact form $\alpha$ used in the proof of \propref{p:D3 V X} are invariant under rotations around the $z$--axis we may assume that the same is true for $V$ and $\alpha(V)$. For $s\in[0,1]$ let 
$$
\beta_s = sr_0\cos(r_0^2\sin^2(\theta))\,d\theta + \sin(r_0^2\sin^2(\theta))\,d\varphi-\alpha(V)\,dt
$$
on $\partial_-R_3=S^2\x S^1$ in spherical coordinates $(\theta,\varphi)$ on $S^2$. According to \eqref{e:perturbed prolong EE}, the contact structure on $\partial_-R_3$ is defined by $\beta_1$ and $\beta_s$ is a contact form for all $s$. By \thmref{t:gray} the contact structures defined by $\beta_0$ and $\beta_1$ are isotopic. The self diffeomorphism $f$ of $S^2\x S^1$ with $f(\theta,\varphi,t)=(\theta,-\varphi,-t)$ preserves $\beta_0$. Moreover $f$ preserves the Legendrian curves $\{p\}\x S^1$ when $p\in\{\varphi=0\}$ and $\alpha(V)$ vanishes at $p$. The orientations are reversed by $f$.

Hence $f$ extends to a diffeomorphism $\widetilde{f}$ of $R_3$ which preserves the contact structure on $\partial_-R_3$ but reverses its orientation and which preserves a Legendrian curve on $\partial_-R_3$ but reverses its orientation. The rotation number of the intersection line field on $\partial_-R_3$ along $\{p\}\x S^1$ is $-|k|$ for the model Engel structures we have obtained so far. If we push forward the model Engel structure we obtained so far using $\widetilde{f}$ we obtain model Engel structures which represent the missing homotopy classes of intersection line fields on $\partial_-R_3$. Thus we have proved the following theorem. 
\begin{thm} \mlabel{t:models on R3}
There are model Engel structures on $R_3$ such that the induced contact structure on $\partial_-R_3$ is $\CC=\ker(\beta_1)$ and all possible orientations of $\CC$ and homotopy classes of intersection line fields are realized.
\end{thm}

\section{Attaching maps for round handles} \mlabel{s:isotopy+choice}
In the following sections we explain how to extend Engel structures from $M$ to $M\cup_\psi R_l$ using a suitable model Engel structure on $R_l$ if $\psi$ is an attaching map for a round handle of index $l=1,2,3$. For this $\psi$ has to have certain properties explained in \secref{s:extend Engel structure} and we will isotope $\psi$ and apply vertical modification of the boundary to $M$ to bring $\psi$ into the desired form. This is discussed in \secref{s:attach 1} to \secref{s:attach 3}. 
\subsection{Extending Engel structures} \label{s:extend Engel structure}
Let $M$ be an Engel manifold with oriented characteristic foliation and transverse boundaries and $R_l$ a round handle of index $l\in\{1,2,3\}$ with a model Engel structure. 
\begin{prop} \mlabel{p:gluing}
Assume that an embedding $\psi : \partial_-R_l \rightarrow \partial_+M$ maps the oriented contact structure $\CC_R$ on $\partial_-R_l$ to the oriented contact structure $\CC_M$ on $\partial_+M$  and preserves intersection line fields. 

Then the Engel structure extends from $M$ to $M'=M\cup_\psi M$ by the model Engel structure on $R_l$. If the Engel structures are oriented, then one obtains an oriented Engel structure on $M'$ if $\psi$ preserves the orientation of the intersection line fields. 
\end{prop}
\begin{proof}
Using the normal form theorem (\thmref{t:normal forms for transversals}) for neighborhoods of transverse hypersurfaces in Engel manifolds we can embed collars $U_R$ of $\partial_-R_l$, respectively $U_M$ of $\partial_+M$, into $\PP\CC_R$, respectively $\PP\CC_M$. We identify $U_R$ and $U_M$ with their image. 

By the second part of \propref{p:Engel diffeos} we obtain an embedding $\widetilde{\psi} : U_R \lra \PP\CC_M$ from $\psi$. Since $\varphi$ preserves intersection line fields, $\widetilde{\psi}$ maps $\partial_-R_l\subset U_R$ to $\partial_+M\subset U_M$ and because $\psi$ preserves the orientation of the contact structure, $\widetilde{\psi}(U_R)$ and $U_M$ lie on opposite sides of the hypersurface $\partial_+M\subset\PP\CC_M$. This proves the claim.
\end{proof}
Strictly speaking $M'$ is a manifold with corners. We smooth these corners by cutting out a piece of $R_l$ such that the boundary of the resulting manifold is still transverse. \propref{p:gluing} can of course be formulated with an Engel manifold instead of the round handle $R_l$.  

The following three sections describe how to choose model Engel structures on $R_1, R_2, R_3$ and how to isotope the attaching map in order to satisfy the conditions in \propref{p:gluing}. It is enough to ensure that the attaching map preserves to homotopy class of the intersection line field since we can then apply vertical modifications of the boundary.

Before we continue let us remark that we will always assume that the attaching map $\psi : \partial_-R_l \lra \partial_+M$ preserves the contact orientations. If this is not the case we can replace $\psi$ by $\psi\circ f$ where $f$ is the diffeomorphism of $\partial_-R_l$ induced by complex conjugation of $S^1$. We obtain diffeomorphic manifolds when we attach $R_l$ to $M$ using $\psi$ or $\psi\circ f$.

\subsection{Attaching maps for round handles of index $1$} \mlabel{s:attach 1}
The model Engel structures on round handles of index $1$ induce contact structures on $\partial_-R_1$ such that the attaching curves $\gamma_\pm=\{\pm 1\}\x\{0\}\x S^1$ are Legendrian. Let $\psi : \partial_-R_1 \lra \partial_+M$ be an embedding preserving the contact orientation. We want to isotope $\psi$ such that the resulting map preserves oriented contact structures and the homotopy class of intersection line fields for a suitable choice of model Engel structure $\DD_{k,m}$. After a vertical modification of $\partial_+M$ we then obtain an Engel structure on $M$ with $R_1$ attached.
 
It is a standard fact in contact topology that every embedded curve in a contact manifold can by isotoped to a Legendrian curve such that the isotopy never leaves an arbitrarily small tubular neighborhood of the original curve, \cite{sympgeo}. Hence we may assume that $\psi(\gamma_\pm)$ are Legendrian curves.

Next we want to isotope $\psi$ such that the resulting map preserves oriented contact structures and the homotopy class of intersection line fields along $\gamma_\pm$ for a suitable choice of model Engel structures on $R_1$. For this we stabilize the curves $\psi(\gamma_\pm)$. We may assume that all stabilizations of $\psi(\gamma_\pm)$ are carried out in disjoint tubular neighborhoods of the original curves. Then one can extend the isotopy of $\psi$ from $\gamma_\pm$ to $\partial_-R_1$ and one obtains a stabilized attaching map.  

Let $\psi_+$ and $\psi_-$ be the restrictions of $\psi$ to $\{x=1\}\x D^1\x S^1$ and $\{x=-1\}\x D^2\x S^1$. By $\fr(\gamma_\pm,m)$ we denote a contact framing of $\gamma_\pm$ if $R_1$ carries the model Engel structure $\DD_{k,m}$.

Since $\psi$ is orientation preserving there exist integers $n_+,n_-$ such that $\psi$ maps a contact framing along $\gamma_\pm$ to a framing representing $n_\pm\cdot\fr(\psi_\pm(\gamma_\pm))$ when $R_1$ carries the model Engel structure $\DD_{k,0}$.
\begin{prop} \mlabel{p:adapt frame and rot}
We can choose a model Engel structure on $R_1$ and stabilize $\psi_\pm$ such that the stabilized maps
\begin{itemize}
\item[(i)] send contact framings of $\gamma_\pm$ to framings of $\psi_\pm(\gamma_\pm)$ which are homotopic to a contact framing, and
\item[(ii)] the rotation numbers of $\DD$ along the stabilized Legendrian curves obtained from $\psi_+(\gamma_+)$ and $\psi_-(\gamma_-)$ are both equal,
\end{itemize}
if and only if 
\begin{equation} \label{e:mod2}
n_+ + \rot(\psi_+(\gamma_+)) \equiv n_- + \rot(\psi_-(\gamma_-)) \mod 2~.
\end{equation}
\end{prop}
\begin{proof}
We equip $R_1$ with the model Engel structure $\DD_{k,0}$. Let $(S,T)$ be a framing of $\gamma)$. Because $\psi$ is orientation preserving, the framings 
$$
m\cdot\big(\psi_{\pm*}(S,T)\big)\hspace{0.5cm} \textrm{ and } \hspace{0.5cm}\psi_{\pm*}\big(m\cdot(S,T)\big)
$$
are homotopic. If we use the model Engel structure $\DD_{k,m}$ instead of $\DD_{k,0}$ on $R_1$ this implies
\begin{align*}
\psi_{\pm*}\big(\fr(\gamma_\pm,m)\big) & = (m+n_\pm)\cdot\fr(\psi_\pm(\gamma_\pm))~.
\end{align*}
In \eqref{e:tb twist} we have determined the effect of positive and negative stabilization on contact framings. Since we want the stabilized embeddings $\widetilde{\psi}_\pm$ to map contact framings of $\gamma_\pm$ to framings of $\widetilde{\psi}_\pm(\gamma_\pm)$ which are homotopic to contact framings, we have to apply positive or negative stabilization $(n_\pm+m)$--times. Hence $n_\pm+m$ has to be non--negative. If $n_+^+, n_+^-,n_-^+,n_-^-\in\N_0$ satisfy
\begin{equation} \label{e:gls framings}
\begin{split}
n_+ + m & = n_+^+ + n_+^-\ge0 \\
n_- + m & = n_-^+ + n_-^-\ge0~,
\end{split}
\end{equation}
then it follows from \eqref{e:rot stab} that the rotation numbers of the stabilized Legendrian curves are given by
\begin{align*}
\rot\left(\left((\sigma_+)^{n_+^+}(\sigma_-)^{n_+^-}\psi_+\right)(\gamma_+)\right) & = \rot(\psi_+(\gamma_+)) + n_+^+ - n_+^- \\
\rot\left(\left((\sigma_+)^{n_-^+}(\sigma_-)^{n_-^-}\psi_-\right)(\gamma_-)\right) & = \rot(\psi_-(\gamma_-)) + n_-^+ - n_-^-~.
\end{align*}
After sufficiently many stabilizations we also want the rotation numbers along the image of $\gamma_+$ and $\gamma_-$ to be equal. This can be achieved if and only if we can solve \eqref{e:gls framings} and
\begin{align} \label{e:gls rot}
n_-^+ - n_-^- - n_+^+ + n_+^- & = \rot(\psi_+(\gamma_+)) - \rot(\psi_-(\gamma_-)) 
\end{align}
with nonnegative integers $n_+^+, n_+^-,n_-^+,n_-^-$ and $m\in\Z$. Then we can take 
\begin{equation*}
k = \rot(\psi_+(\gamma_+)) + n_+^+ - n_+^- = \rot(\psi_-(\gamma_-)) + n_-^+ - n_-^-~. 
\end{equation*}
From \eqref{e:gls framings} and \eqref{e:gls rot} we obtain the condition \eqref{e:mod2}. Therefore this condition is necessary. 

Conversely, if \eqref{e:mod2} is satisfied, then the under--determined system of equations \eqref{e:rot stab} and \eqref{e:gls rot} admits solutions in $\Z$. If we choose $m$ large enough we can achieve $n_+^+,n_+^-,n_-^+,n_-^-\in\N_0$.
\end{proof}  
Note that in this proof we used only stabilization to modify Legendrian curves. Attaching a bypass to a convex surface along a segment of a Legendrian knot contained in the surface can be viewed as an inverse procedure of stabilization. Under the assumption that the contact structure on $\partial_+M$ is overtwisted it is possible to prove \propref{p:adapt frame and rot} for model Engel structures with a fixed contact framing, i.e. without choosing $m$ big enough as we did above.  

Let us explain the meaning of \eqref{e:mod2} in more topological terms. Consider an orientation preserving attaching map $\psi : \partial_-R_1 \rightarrow \partial_+M$. On $TM$ we have the Engel framing and we may assume that its components except of the component spanning $\WW$ form a framing of the tangent bundle of $\partial_+M$. Using $\psi$ we pull back this trivialization to $\partial_-R_1$ and add a vector field which points inwards and is transverse to $\partial_-R_1$. We obtain a framing of $TR_1$ along $\partial_-R_1$.

It is possible to isotope $\psi$ such that we can extend the Engel structure from $M$ to $M$ with $R_1$ attached only if we can extend the Engel framing from $M$ to $M\cup_\psi R_1$. This condition depends only on the Engel framing and the isotopy class of $\psi$. If $\psi(\gamma_\pm)$ are Legendrian one can show that such an extension exists if and only if \eqref{e:mod2} is satisfied. 
 
On the other hand extending a framing from $M$ to $M\cup_\psi R_1$ is a purely topological problem. Thus if we start with an attaching map $\psi_\pm$ and end up with a map violating \eqref{e:mod2} then it is not possible to construct an Engel structure on $M\cup_\psi R_1$ such that the Engel framing on $M$ is homotopic to the Engel framing induced by $\DD$. 

Now if $\psi : \partial_-R_1 \lra \partial_+M$ is an embedding such that $\psi(\gamma_\pm)$ are Legendrian curves and $\varphi$ preserves contact framings and rotation numbers along $\gamma_\pm$ then by \lemref{l:tub nbhd of leg curves} we can isotope $\psi$ relative to $\gamma_\pm$ such that the resulting map preserves the contact structure on a tubular neighborhood of $\gamma_\pm$. Using \remref{r:contact contraction} we can isotope the new attaching map such that its image lies in this tubular neighborhood. 

\begin{thm}  \mlabel{t:construct attaching maps}
Let $M$ be an oriented Engel manifold with transverse boundary and assume that $\psi : \partial_-R_1 \lra \partial_+M$ is an embedding and the trivialization of $TM$ induced by the Engel structure can be extended to $M\cup_\varphi R_1$.

Then there is a model Engel structure on $R_1$ such that $\psi$ is isotopic to an embedding $\widetilde{\psi}$ which preserves oriented contact structures and the homotopy class of the intersection line field.
\end{thm}

%%%%%%%%%%%%%%%%%%%%%%%%%%%%%%%%%%%%%%%%%%%%%%%%%%%%%%%%%%%%
\subsection{Attaching maps for round handles of index $2$} \mlabel{s:flex tori}

All our model Engel structures on round $2$--handles induce the same singular foliation on $T^2_0=\partial D^2\x\{0\}\x S^1\subset\partial_-R_2$. Now suppose that $M$ is an Engel manifold with transverse boundary and $\psi : \partial_-R_2 \lra \partial_+M$ is an attaching map which preserves the orientations induced by the contact structures.

If we want to attach $R_2$ to $M$ and extend the Engel structure from $M$ to $M\cup R_2$, we have to ensure that the attaching map preserves contact structures. By \thmref{t:divide impera} and \thmref{t:giroux unique} together with \remref{r:contact contraction 2} it suffices to modify $\psi$ such that after the deformation, the image of $T_0^2$ is convex and the attaching map preserves singular foliation. Recall that the dividing set of $T^2_0$ consists of two homotopically non--trivial circles.   

Let $N=\partial_+M$ and $T^2=\psi(T^2_0)$. We assume that the contact structure $\CC$ on $N$ is overtwisted. Using \lemref{l:attach} and \propref{p:bypass existence} we can bring $T^2$ into the desired form. It is clear how to obtain the desired isotopy of $\psi$ from this.
\begin{thm} \mlabel{t:flexible tori}
Let $T^2$ be an embedded torus in an overtwisted contact manifold $(N,\CC)$. Assume that $\CC$ is orientable and that the Euler class of the restriction of $\CC$ to $T^2$ is zero. Then we can isotope $T^2$ such that after the isotopy the singular foliation on the torus is in standard form. Moreover we can prescribe the slope of the dividing curves.

After the isotopy, the complement of a tubular neighborhood of $T^2$ contains an overtwisted disc. 
\end{thm} 
\begin{proof}
It suffices to find a convex torus which is isotopic to the original one such that the dividing set consists of two homotopically non--trivial components which have the desired slope. Using the Giroux flexibility theorem (\thmref{t:giroux flex}) one can arrange the singular foliation on $T^2$ such that it has standard form. We will frequently use \propref{p:bypass existence}. The following figures represent the dividing set on a torus before the bypass attachment. The thickened arc represents the attaching curve $\gamma_1$ of the bypass.

$1^{st}$ Step: 
Let $D_{ot}$ be a convex overtwisted disc. We perturb the embedding of $T^2$ such that it becomes transverse to $D_{ot}$. Using an extension of a radial vector field on $D_{ot}$ we can isotope $T^2$ such that after the isotopy $T^2\cap D_{ot}=\emptyset$ and the resulting torus is convex.

Since $D_{ot}$ is convex there is a neighborhood which is foliated by overtwisted discs. In the following we will always ensure that after each modification of the embedding of $T^2$, there is an overtwisted disc disjoint from the deformed torus: If $D$ is a bypass for $T^2$, we choose the neighborhood of $T^2\cup D$ in \lemref{l:attach} so small that its complement still contains overtwisted discs. 

$2^{nd}$ Step:
In this step we remove all homotopically non--trivial components of the dividing set. If there are no such components we continue with step 3.

If the dividing set contains more than two homotopically non--trivial components, then we reduce the number of its components of the dividing set using the bypass attachments in the left part of \figref{b:weg} often enough. 
\begin{figure}[htb]
\begin{center}
\includegraphics{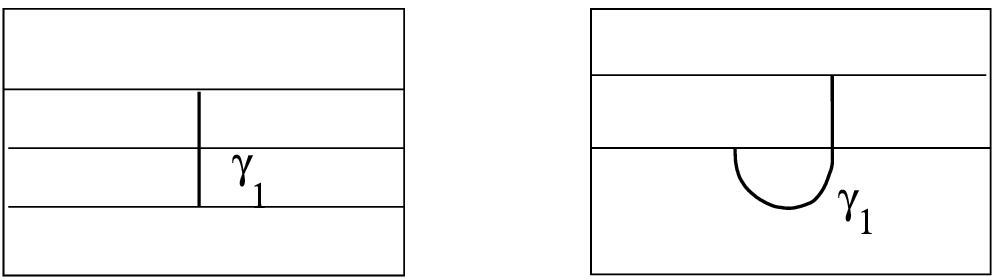}
\caption{\label{b:weg}}
\end{center}
\end{figure}
We end up with a dividing set which contains two homotopically non--trivial curves. We remove these components with the bypass attachment in the right part of \figref{b:weg}

$3^{rd}$ Step:
Using the bypass attachment in \figref{b:neuecomp}, we obtain two new components of the dividing set. Their slope depends on the attaching curve of the bypass. For a given identification $T^2\simeq S^1\x S^1$, we can achieve that the new components of the dividing set are isotopic to $\{p\}\x S^1$ for $p\in S^1$. The dashed curve represents this circle. 
\begin{figure}[htb] 
\begin{center}
\includegraphics{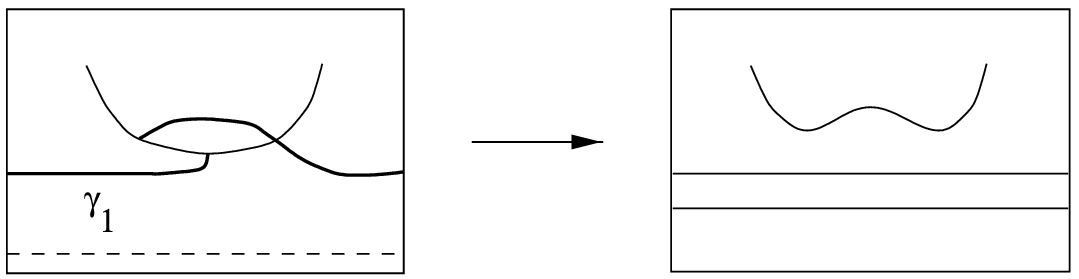}
\end{center}
\caption{\label{b:neuecomp}}
\end{figure}

$4^{th}$ Step:
We are left with a convex torus whose dividing set contains exactly two homotopically non--trivial dividing curves $\sigma_1,\sigma_2$ with the desired slope. If this is the entire dividing set we are done. Otherwise we consider the two annuli $T^2\setminus(\sigma_1\cup\sigma_2)$. 

We claim that if only one of these annuli contains other components of the dividing set $\Gamma$, then there is at least one component of $\Gamma$ which bounds a disc $\widetilde{D}$ containing another component of $\Gamma$. Assume that this is not true. Then $T^2\setminus\Gamma$ contains $r>0$ discs, one annulus and one annulus with $r$ holes. The Euler number $\chi(\CC,T^2)$ of the restriction of $\CC$ to $T^2$ is
\begin{equation} \label{e:eulerchar fuer einf torus}
\chi(\CC,T^2) = \chi(T^2_+)-\chi(T^2_-) = \pm 2r\neq 0~.
\end{equation}
by \eqref{e:eulerchar}. The sign depends on the orientations of $T^2$ and of the contact structure. But \eqref{e:eulerchar fuer einf torus} contradicts our assumption on the Euler class of $\CC$. In order to reduce the number of connected components of $\Gamma$ we perform a bypass attachment as indicated in the left part of \figref{b:otweg}. Notice that this does not affect the homotopically non--trivial dividing curves.
\begin{figure}[htb]
\begin{center}
\includegraphics{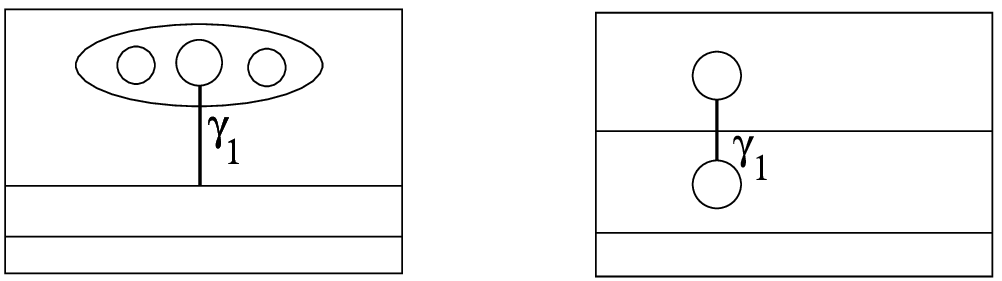}
\caption{\label{b:otweg}}
\end{center}
\end{figure}

If both annuli $T^2\setminus(\sigma_1\cup\sigma_2)$ contain connected components of $\Gamma$ we reduce the number of components using the bypass attachment in the right part of \figref{b:otweg}. If do this often enough we end up with the desired configuration of dividing curves on $T^2$. 
\end{proof}
A bypass attachment also affects framings. For our purpose, it is enough to show this for a particular bypass attachment. K.~Honda described this effect in more detail, cf. Proposition 4.7 in \cite{honda} where one can find a proof of the following lemma.
\begin{figure}[htb] 
\begin{center}
\includegraphics{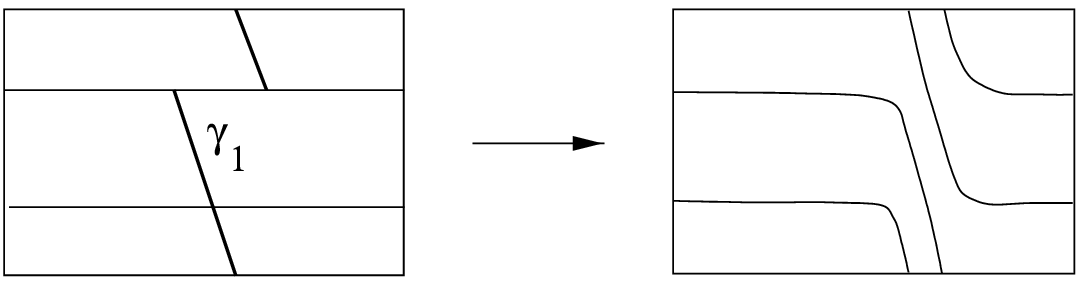}
\end{center}
\caption{\label{b:dehn}}
\end{figure}

\begin{lem} \mlabel{l:dehn}
Let $T^2=S^1\x S^1$ be a torus in standard form contained in a contact manifold such that the Legendrian divides are isotopic to $\{1\}\x S^1$. Let $X$ be a nowhere vanishing Legendrian vector field on the contact manifold. 

Assume that the rotation number of $X$ along the Legendrian divides is zero and that it is even along $S^1\x\{0\}$. We attach a bypass as in \figref{b:dehn} to $T^2$ and bring the characteristic foliation in standard form such that the Legendrian ruling is still tangent to the foliation from the first factor in $S^1\x S^1$. 

Then the rotation number along the Legendrian divides in the isotoped torus is odd (and therefore non--zero) while the rotation number along the Legendrian curve $S^1\x\{0\}$ remains even.  
\end{lem}
 
\begin{prop} \mlabel{p:attaching round 2-handles}
Assume that $M$ carries an oriented Engel structure with oriented characteristic foliation such that the boundary of $M$ is transverse. Moreover, assume that the contact structure on $\partial_+M$ is overtwisted. Let $\psi : \partial_-R_2 \lra \partial_+M$ be an attaching map.

We can extend the Engel structure from $M$ to $M\cup R_2$ using a model Engel structure if and only if the Engel framing on $M$ extends to a trivialization of $M\cup R_2$ over the $2$--cell $e_2=D^2\x\{0\}\x\{1\}\subset R_2$. The Engel structure on $M\cup R_2$ can be chosen such that the contact structure on $\partial_+(M\cup R_2)$ again overtwisted.
\end{prop}
\begin{proof} 
The $2$--cell $e_2$ is of course attached using the restriction of $\psi$ to $\partial e_2$. Clearly, if the Engel trivialization can not be extended over $e_2$, then it is impossible to extend the Engel structure to $M\cup_{\psi_2} R_2$. 

The contact structure $\CC$ on $\partial_+M$ is oriented and the intersection line field yields a nowhere vanishing Legendrian vector field. Therefore the Euler characteristic of $\CC$ viewed as a bundle is zero. We equip $R_2$ with one of the model Engel structures from \secref{s:modelle index 2}. By \thmref{t:flexible tori} and \remref{r:contact contraction 2} we can assume that $\psi$ preserves contact structures. In particular, $\psi(\partial e_2)$ is Legendrian. 

We first show that under the assumption of the proposition, the rotation number along $\psi(\partial e_2)$ is even. For this we homotop the Engel trivialization on $M$ such that the only component $W$ of the framing which is not tangent to $\partial_+M$ spans the characteristic foliation. In this way we obtain a framing of $\partial_+M$. Using the attaching map to pull back this framing, we obtain a framing on $\partial_-R_2$. If we add a vector field along $\partial_-R_2$ which spans and orients the characteristic foliation we obtain a framing of $TR_2$ along $\partial_-R_2$. 

Since $\psi$ preserves contact structures, the pull back of the component of the Engel framing which is orthogonal to the contact structure on $\partial_+M$ is transverse to the contact structure on $\partial_-R_2$. Without loss of generality, we assume that these components of the Engel framings are preserved by the attaching map. Then the pullback framing and the Engel framing on $R_2$ have two components in common. When we want to compare the pull back framing with the Engel trivialization along $e_2$ it is therefore enough to consider the rotation numbers along $\partial e_2$. 

By (iii) of \propref{p:zwei henkel modell}, the rotation number along $\partial e_2$ induced by the model Engel structure on $R_2$ is even. If the rotation number of the pull back framing along $\partial e_2$ is odd, then the pull back framing and the Engel framing are not homotopic along $\partial e_2$ and the pull back framing can not be extended over the disc $D^2\x\{0\}\subset h_2$. This contradicts the assumption. Hence the rotation number along $\psi(\partial e_2)$ is even.  

If the rotation number along the Legendrian divides of $\psi(T_0^2)$ is not zero, then by \propref{p:zwei henkel modell} we can choose a suitable model Engel structure such that the isotoped attaching map preserves the homotopy class of the intersection line fields. 

If the rotation number along the Legendrian divides of $\psi(T_0^2)$ is zero, then we attach a bypass as in \figref{b:dehn} to isotope $\psi$ once again. By \lemref{l:dehn}, the rotation number along the Legendrian rulings of the image of $T^0_0$ remains even and the rotation number along the Legendrian divides is non--zero. 

The new attaching map preserves the singular foliation on $T_0^2$ after one pushes forward the model Engel structures from \propref{p:zwei henkel modell} with the diffeomorphism $\Theta_1$ from \eqref{e:Theta}.

When we restrict this diffeomorphism to $\partial D^2\x\{0\}\x S^1$ then we obtain a right handed Dehn twist. Now we can choose the model Engel structure such that the new attaching map preserves oriented contact structures and the homotopy class of the intersection line field. After a vertical modification of the boundary of $M$, we can extend the Engel structure over $R_2$ using a model Engel structure from \secref{s:modelle index 2}.
\end{proof}

%%%%%%%%%%%%%%%%%%%%%%%%%%%%%%%%%%%%%%%%%%%%%%%%%%%%%%%%%%%%
\subsection{Attaching maps for round handles of index $3$} \mlabel{s:attach 3}

Let $\psi : \partial_-R_3 \lra \partial_+M$ be an orientation preserving attaching map for a round $3$--handle carrying a model Engel structure from \secref{ss:model 3-handles}. The model Engel structures on $R_3$ induce overtwisted contact structures on $\partial_-R_3$. By \thmref{t:eliashbergs wahnsinn} we can isotope $\psi$ to a contact diffeomorphism if and only if $\psi$ is orientation preserving and the image of the contact structure on $\partial_-R_3$ is homotopic to the contact structure on $\partial_+M$ as a plane field. 

The contact structure on $\partial_-R_3$ is orientable and it extends to $R_3$ as a plane field. The following lemma shows that this determines the homotopy class of the plane field completely.
\begin{lem} \mlabel{l:unique homotopy class of plane fields}
There is a unique homotopy class of orientable plane fields on $S^2\x S^1=\partial D^3\x S^1$ which extends to $D^3\x S^1$. 
\end{lem}
\begin{proof} 
Recall from \cite{hirzhopf} that the Grassmannian of oriented planes in $\R^3$, respectively $\R^4$, is $\textrm{Gr}_2(3)\simeq S^2$, respectively $\textrm{Gr}_2(4)\simeq S^2\x S^2$ and that the inclusion $\R^3\lra\R^4$ induces the diagonal map
$$
\Delta : \textrm{Gr}_2(3)\simeq S^2 \lra S^2\x S^2\simeq \textrm{Gr}_2(4)~.
$$
Let $\CC_0$ and $\CC_1$ be two plane fields on $S^2\x S^1$ that extend to the interior of $D^3\x S^1$. We view $\CC_0,\CC_1$ as maps from $S^2\x S^1$ to $\textrm{Gr}_2(3)$ and their extensions as maps from $D^3\x S^1$ to $\textrm{Gr}_2(4)$. 

Because $\{0\}\x S^1$ is a strong deformation retract of $D^3\x S^1$ and $\textrm{Gr}_2(4)$ is simply connected, there is a homotopy between the extensions of $\CC_0$ and $\CC_1$ through plane fields in $T(D^3\x S^1)$. Using the projection of $\textrm{Gr}_2(4)\simeq S^2\x S^2$ onto the first factor, we obtain a homotopy between $\CC_0$ and $\CC_1$ in $T(S^2\x S^1)$.
\end{proof}
After we have isotoped $\psi$ to a contact diffeomorphism we can choose a model Engel structure on $R_3$ such that the orientation of the contact structure and the homotopy class of the intersection line field are preserved by the isotoped attaching map, cf.~\thmref{t:models on R3}.
\begin{prop} \mlabel{p:extending over R3}
The Engel structure extends from $M$ to $M\cup_\psi R_3$ if and only if the contact structure on $\partial_+M$ extends to $M\cup_\psi R_3$ as a plane field.
\end{prop} 

%%%%%%%%%%%%%%%%%%%%%%%%%%%%%%%%%%%%%%%%%%%%%%%%%%%%%%%%%%%%
\section{Existence theorems} \mlabel{s:proofs}

Now we combine the tools from the previous sections to prove existence theorems for Engel structures. In \secref{s:end proof} we show that every parallelizable manifold admits an orientable Engel structure. In \secref{s:conn sum} we explain the construction of Engel structures on the connected sum of two Engel manifolds with $S^2\x S^2$ if the characteristic foliations of the original Engel structures satisfy certain conditions. 

%%%%%%%%%%%%%%%%%%%%%%%%%%%%%%%%%%%%%%%%%%%%%%%%%%%%%%%%%%%%
\subsection{Manifolds with trivial tangent bundle} \mlabel{s:end proof}

In the proof of our main result \thmref{T:EXIST} we use the model Engel structures from \secref{s:round model} and the results from \secref{s:isotopy+choice} showing that we can extend an Engel structure from $M'\subset M$ to $M'\cup R_l$ if $M'$ has transverse boundary and the Engel framing extends to a framing of $M'\cup R_l$. The remaining problem is to ensure that after we have constructed an Engel structure on $M'\cup R_l$, the new Engel framing extends from $M'\cup R_l$ to the whole of $M$. This allows us to perform the handle attachments successively. 
\begin{thm} \mlabel{T:EXIST}
Every parallelizable $4$--manifold admits an orientable Engel structure. 
\end{thm}
\begin{proof}
Let $M$ be a closed parallelizable manifold of dimension $4$ and fix a trivialization $TM\simeq M\x\R^4$ of $TM$. We consider a round handle decomposition 
$$
M =\left(\ldots\left(\left(\ldots\left( R_0\cup_{\psi_1^1} R_1^1\right) \ldots\cup_{\psi_1^{r_1}}R_1^{r_1}\right) \cup_{\psi_2^1}R_2^1\right)\ldots\cup_{\psi_2^{r_2}}R_2^{r_2}\right)\cup_{\psi_3} R_3
$$
of $M$ such that there is exactly one round $3$--handle and one round $0$--handle. Such a decomposition exists by \thmref{t:round decomp euler 0}. We will frequently isotope the attaching maps but this will not be reflected in the notation.

We start with the round handle of index $0$. By \propref{p:index 0 models} we can choose a model Engel structure on $R_0$ such that the Engel framing on $R_0$ is homotopic to the original framing. In particular the Engel framing extends from $R_0$ to $M$. The contact structure on $\partial_+R_0$ is overtwisted by construction.

Let $M_1^{i-1}$ be the round handle body obtained from $R_0, R_1^{1}, \ldots, R_1^{i-1}$. Assume that we have constructed an Engel structure on $M_1^{i-1}$ such that the contact structure on $\partial_+M_1^{i-1}$ is overtwisted. Assume moreover that throughout this process we have homotoped the trivialization of $M$ we chose at the beginning such that it coincides with the Engel trivialization on $M_1^{i-1}$.

Then the Engel trivialization on $M_1^{i-1}$ can be extended to $M_1^{i-1}\cup_{\psi_1^i} R_1^i$. By \thmref{t:construct attaching maps}, we can isotope $\psi_1^i$ to an attaching map such that the Engel structure on $M_1^{i-1}$ extends to an Engel structure on $M_1^i=M_1^i\cup_{\psi_1^i} R_1^i$ using a model Engel structure on $R_1$ from \secref{s:index one models}. In order to ensure that the contact structure on $\partial_+M_1^{i}$ is again overtwisted, we isotope $\psi_1^i$ before the application of \thmref{t:construct attaching maps} such that its image is disjoint from an overtwisted disc in $\partial_+M_1^{i-1}$. 

Let $\gamma_\pm=\{\pm 1\}\x\{0\}\x S^1$ be the attaching curves of $R_1^i$. Assume that $\psi_1^i(\gamma_\pm)$ is transverse to an overtwisted disc $D_{ot}$ and let $p$ be a point on $D_{ot}$ which does not lie on $\psi_1^i(\gamma_\pm)$. Then use the flow of a radial vector field centered at $p$ to isotope $\psi_1^i$ such that the image of $\gamma_\pm$ becomes disjoint from $D_{ot}$. The remaining steps in the modification of $\psi$, like making the attaching curves Legendrian and performing stabilizations, can be carried out in a small tubular neighborhood which is also disjoint from $D_{ot}$.

Next we compare the Engel trivialization and the original trivialization of $M$ on $M_1^i$ relative to $M_1^{i-1}$. The cylinder 
$$
D^1\x \{0\}\x S^1\subset R_1^i= D^1\x D^2\x S^1
$$ 
can be decomposed into a $1$--cell $e_1=D^1\x\{0\}\x\{1\}$ and a $2$--cell $e_2$. The $1$--cell is attached to $M^{i-1}$ using the restriction of $\psi_1^i$ to $\partial e_1$ and $e_2$ is attached to $M_1^{i-1}\cup e_1$. Since $\pi_1(\SO(4))=\Z_2$ there are two homotopy classes of oriented framings of $TM$ along $e_1$.

If necessary, we modify the model Engel structure on $R_1^i$ such that the new Engel framing is homotopic to the given trivialization along $e_1$ relative to the endpoints of $e_1$. For this let $\rho : D^1 \lra [0,2\pi]$ be a smooth function which is constant near the boundary, $\rho(-1)=0$ and $\rho(1)=2\pi$. If we push forward the model Engel structure on $R_1$ using the diffeomorphism 
\begin{align*}
F_1 : R_1^i=D^1\x D^2\x S^1 & \lra D^1\x D^2\x S^1=R_1^i \\
          (x,y_1,y_2,t) & \lmt \begin{array}{c}(x,\cos(\rho(x))y_1+\sin(\rho(x))y_2, \\ -\sin(\rho(x))y_1 + \cos(\rho(x))y_2, t) \end{array}~,
\end{align*}
then we still obtain a smooth Engel structure on $M_1^i$ by the choice of $\rho$. Because $F_1$ interchanges the two homotopy classes of framings along $e_1$ relative $\partial e_1$, the trivialization induced by the new Engel structure on $R_1$ is now homotopic to the given trivialization of $TM$ along $e_1$ relative $\partial e_1$. We homotop the original framing such it coincides with the Engel framing on $M_1^{i-1}\cup e_1$. Since $\pi_2(\SO(4))$ is trivial the same is true for $M^{i-1}_1\cup e_1\cup e_2$.  
 
Thus the Engel framing on $M_1^{i}$ is now homotopic to the original framing relative to $M_1^{i-1}$. In particular it extends from $M_1^{i}$ to $M$. Hence we can iterate the procedure to attach round handles of index $1$. During this process we ensure that the contact structure on $\partial_+M_1^{i}$ is overtwisted for all $i$ by choosing the attaching regions disjoint from $D_{ot}$. 

In the next step we attach round $2$--handles. Assume that we have already attached the first $i-1$ round $2$--handles such that on the resulting handle body $M_2^{i-1}$ there is an Engel structure extending the Engel structure on $M_1$. In addition we assume that the contact structure on $\partial_+M_2^{i-1}$ is overtwisted. Consider the attaching map 
$$ 
\psi_2^i : \partial_-R_2^i \lra \partial_+M_2^{i-1}~.
$$
of $R_2^i$. The contact structure on $\partial_+M_2^{i-1}$ is oriented and it has an orientable section, namely the intersection line field. Thus the Euler class of the contact structure, viewed as a bundle, vanishes. By assumption, the contact structure is overtwisted. According to \thmref{t:flexible tori} and \thmref{t:giroux flex} we can isotope $\psi_2^i$ such that the resulting map preserves singular foliations. Recall that in \thmref{t:flexible tori} we also ensured that the attaching region of $R_2^i$ is contained in a neighborhood of $\psi_2^i(T^2_0)$ which is disjoint from some overtwisted disc.

In order to find a model Engel structure on $R_2^i$ which extends the Engel structure on $M_2^{i-1}$ to an Engel structure on $M_2^i=M_2^{i-1}\cup_{\psi_2^i}R_2^i$ we have to prove the following claim. 

{\it Claim 1 : The Engel trivialization on $M_2^{i-1}$ extends to a trivialization of $TM$ over $e_2\subset R_2^i$. The extension is unique up to homotopy relative $\partial e_2$.}
\begin{proof}[Proof of Claim 1]
First we decompose each round $2$--handle $R_2^j, j\le i$ into one ordinary handle $h_2^j$ of index $2$ and one ordinary handle $h_3^j$ of index $3$ as in \eqref{e:decomp Rk}. Then we rearrange the handles $h_2^j, h_3^j, j\le i$ such that the ordinary handles of index $2$ are attached to $M_1$ independently and the remaining ordinary $3$--handles are attached to the resulting handle body.

The Engel trivialization on $M_2^{i-1}$ extends over $h_2^{i}$ {\sl before} we rearrange the handle decomposition if and only if the Engel trivialization extends from $M_1$ over $h_2^{i}$ {\sl after} the rearrangement. By construction the Engel framing on $M_1$ is homotopic to the framing of $M$ we chose at the beginning of the proof. In particular the Engel framing extends from $M_1$ to $h_2^{i}$ after the rearrangement of the handles $h_2^j, h_2^j, j\le i$. 

This proves that the Engel framing can be extended from $M_2^{i-1}$ to $M_2^{i-1}\cup e_2^i$. This extension is unique up to homotopy relative to $\partial e_2^i$ since $\pi_2(\SO(4))=\{0\}$.
\end{proof}
By \propref{p:attaching round 2-handles} we can extend the Engel structure from $M_2^{i-1}$ to $M_2^i$ using a model Engel structure from \secref{s:modelle index 2} on $R_2^i$. Since the contact structure on $\partial_+M^i_2$ is again overtwisted we can iterate the procedure. When the last round $2$--handle is attached we have constructed an Engel structure on $M_2$. In order to finish the construction we have to extend the Engel structure over the round $3$--handle $R_3$. 

{\it Claim 2 : The Engel trivialization extends from $M_2$ to $M$.}
\begin{proof}[Proof of Claim 2]
We decompose all round $2$--handles into ordinary handles $h_2^j,h_3^j$ for $1\le j\le r_2$ of index $2$ and $3$ and we rearrange the handles such that the $2$--handles are attached to $M_1$. We also decompose $R_3$ into an ordinary $3$--handle $\widehat{h}_3$ and one ordinary $4$--handle $\widehat{h}_4$. In Claim 1 we have shown that the Engel trivialization on $M_1$ extends to $M_1\cup h_2^1\cup\ldots\cup h_2^{r_2}$ and that all such extensions are homotopic. Therefore, the Engel trivialization on $M_1\cup h_2^1\cup\ldots h_2^{r_2}$ also extends to $M$.

Next we reduce the problem to trivializations of bundles of rank $3$. The first component $W$ of the Engel trivialization is transverse to $\partial M_2$ by construction. Thus $W$ extends to a vector field without zeroes on $M$. We equip $M$ with an almost quaternionic structure such that the Engel framing and $W,IW,JW,KW$ coincide on $M_1\cup h_2^1\ldots\cup h_2^{r_2}$. Then we can choose a trivialization of the orthogonal complement $\WW^{\bot}$ of $W$ in $M$. (This trick can be found in \cite{geiges}.) For the proof of Claim 2 it suffices to show that we can extend the trivialization of $\WW^{\bot}$ induced by the Engel structure.

On the $2$--skeleton of $M$ the $\SO(3)$--bundle $\WW^{\bot}$ is trivial. Therefore we can lift it to an $\textrm{SU}(2)$--bundle. (Recall that $\textrm{Spin}(3)=\textrm{SU}(2)=S^3$.) Since $\pi_2(\textrm{SU}(2))$ is trivial, the trivialization of $\WW^{\bot}$ induced by the Engel structure extends from $M_2$ to $\widehat{h}_3$. We fix such an extension. 

The obstruction for the extension of the trivialization from $M_2\cup\widehat{h}_3$ to $M$ of $\WW^{\bot}$ induced by the Engel structure is a cocycle $x$ in the cellular cochain group $C^4(M,M_2\cup\widehat{h}_3;\pi_3(\textrm{SU}(2))=\Z)$ and $x$ may depend on the choice of extensions of the trivialization over the $3$--handles. The cocycle $x$ represents a class $[x]\in H^4(M,M_2\cup\widehat{h}_3;\Z)$ which does not depend on the choice of trivializations on the $3$--handles because it is a primary obstruction, cf. Theorem 34.2 in \cite{steenrod}.

We have already shown that the Engel trivialization extends from the $2$--skeleton to the whole of $M$. Thus $[x]=0$. Since there is exactly one $4$--handle we have $C^4(M,M_2\cup\widehat{h}_3;\Z)=H^4(M,M_2\cup\widehat{h}_3;\Z)$ and hence $x=0$. This implies that the Engel trivialization on $M_2$ extends to $M$. 
\end{proof}
Since the Engel framing extends from $M_2$ to $M$ the same is true for the contact structure on $\partial_+M_2$ viewed as a plane field. The contact structure on $\partial_+M_2$ is overtwisted by construction. By \propref{p:extending over R3} we can extend the Engel structure from $M_2$ to $M$. This finishes the proof of the theorem. 
\end{proof}
Let us compare the proof of \thmref{T:EXIST} with the following characterization of parallelizable $4$--manifolds.
\begin{thm}[Hirzebruch, Hopf, \cite{hirzhopf}]\mlabel{t:trivial charkl}
An orientable $4$--manifold has trivial tangent bundle if and only if
\begin{itemize}
\item[(i)] the Euler characteristic vanishes,
\item[(ii)] the second Stiefel--Whitney class $w_2(M)$ is zero and
\item[(iii)] the signature of $M$ is zero.
\end{itemize}
\end{thm}
Our construction relies on round handle decompositions and therefore condition (i) in \thmref{t:round decomp euler 0} is used at all stages of the proof. The second Stiefel--Whitney class $w_2(M)$ of an orientable $4$--manifold $M$ is zero if and only if $TM$ is trivial on the $2$--skeleton of $M$. Thus condition (ii) of \thmref{t:trivial charkl} is used when we attach round handles of index $1$ and $2$ in order to establish certain properties of rotation numbers. It is also used when we lift $\WW^{\bot}$ to an $\textrm{SU}(2)$--bundle in Claim 2. 

Finally the obstruction for the extension of a section of the $\textrm{SU}(2)$--bundle appearing at the end of the proof of \thmref{T:EXIST} can be viewed as the primary obstruction to the construction of a section of the $\textrm{SU}(2)$--bundle over $M$ since all sections of a $\textrm{SU}(2)$--bundle over $M_1\cup h_2^1\ldots\cup h_2^{r_2}$ are homotopic. According to \cite{gompf} p.~31, an $\textrm{SU(2)}$--bundle is trivial if and only if its second Chern class vanishes. On the other hand the second Chern class of the $\textrm{SU}(2)$--bundle is $-p_1(\WW^\bot)/4$. Since $p_1(TM)=p_1(\WW^\bot)$, the vanishing of $x$ corresponds to (iii) by the signature theorem of Hirzebruch.

\subsection{Connected sums of Engel manifolds} \mlabel{s:conn sum}

Let $M_1, M_2$ be manifolds with Engel structures $\DD_1,\DD_2$. The connected sum $M_1\# M_2$ does not admit an Engel structure because its Euler characteristic is $-2$. On the other hand the Euler characteristic of $M_1\# M_2\# (S^2\x S^2)$ is zero. Moreover, if $M_1, M_2$ have trivial tangent bundle, then the same is true for $M_1\# M_2 \# (S^2\x S^2)$ by \thmref{t:trivial charkl}. 

Therefore it is natural to try to construct an Engel structure on $M_1\# M_2 \# (S^2\x S^2)$ from $\DD_1$ and $\DD_2$ such that the resulting Engel structure coincides with the original Engel structures one large parts of $M_1$ and $M_2$. In this section we show that this is possible under certain assumptions on $\DD_1$ and $\DD_2$.

For this we use model Engel structures on round handles of index $1$ and $2$. As model Engel structure $\DD^{(1)}$ on $R_1$ we use $\DD_1$ from \propref{p:1 modelle}. On $R_2$ we will use a model Engel structure $\DD^{(2)}$ which did not appear yet. We define $\DD^{(2)}$ as the span of 
\begin{align*}
W & = \rve{t}-1/2\,y_1\rve{y_1}-y_2\rve{y_2}+1/2\,x\rve{x} \\
X & = \cos(t)\left(y_2\rve{y_1}+\rve{x}\right) + \sin(t)\left(1/2\,x\rve{y_1} + \rve{y_2}\right)~.
\end{align*}
This model Engel structure can be obtained using \propref{p:perturbed prolong}. Its characteristic foliation is spanned and oriented by $W$. 

The model Engel structures $\DD^{(1)}$ and $\DD^{(2)}$ are very similar when one identifies $R_1$ and $R_2$ in the obvious way. The curves $\gamma_\pm=\{0\}\x\{\pm 1\}\x S^1$ in $\partial_+R_2\simeq\partial_-R_1$ are Legendrian and the rotation number of the intersection line field along these curves is $-1$ for both $\DD^{(1)}$ on $R_1$ and $\DD^{(2)}$ on $R_2$. On a neighborhood of the boundary of $R_1=R_2$ the even contact structures induced by $\DD^{(1)}$ and $\DD^{(2)}$ are homotopic through the family of even contact structures defined by 
\begin{equation} \label{e:alpha s}
\beta_s = -(1-2s)dy_1+1/2\,y_1\, dt-y_2\,dx-1/2\,x\,dy_2
\end{equation}
with $s\in[0,1]$. By \eqref{e:runde 1-henkel} the even contact structure on $R_1$ is defined by $\beta_0$ while the even contact structure on $R_2$ is defined by $\beta_1$. The characteristic foliations of the even contact structures defined by $\beta_s$ are transverse to $\partial_-R_1\simeq\partial_+R_2$ and $\partial_+R_1\simeq\partial_-R_2$ for all $s$. Hence $\beta_s$ induces a family of contact forms on $\partial_-R_1$ and $\partial_+R_1$. 

We can also obtain contact embeddings of $\partial_+R_2$ from attaching maps of $R_1$. Consider an Engel manifold $M$ with transverse boundary and let $\psi_1 : \partial_-R_1 \lra \partial_+M$ be an embedding which preserves oriented contact structures. By \remref{r:contact contraction} we can isotope $\psi_1$ such that the image of the resulting embedding lies in a small neighbourhood of $\psi_1(\gamma_\pm)$. The new embedding is again called $\psi_1$. For a suitable cut off function $\rho$ 
\begin{equation} \label{e:alpha rho}
\widetilde{\beta}_s = -(1-2\rho s)dy_1+1/2\,y_1\, dt-y_2\,dx-1/2\,x\,dy_2
\end{equation}
defines a family of contact structures on $M$ which is constant away from a neighborhood of $\psi_1(\gamma_\pm)$ and which is defined by $\beta_s$ near $\psi_1(\gamma_\pm)$. This homotopy through contact structures has compact support. Using Gray's theorem (\thmref{t:gray}) we obtain an isotopy $f_s$ of $M$ such that $f_1^*\widetilde{\beta}_1$ is a positive multiple of $\widetilde{\beta}_0$. Moreover, $f_s$ maps the Legendrian curves $\psi_1(\gamma_\pm)$ to themselves for all $s$.

The orientation of the contact structure induced by $\DD^{(1)}$ on $\partial_-R_1$, respectively $\DD^{(2)}$ on $\partial_+R_2$, is given by the restriction of $d\beta_0$, respectively $d\beta_1$, to the contact structure. This implies that 
$$
\psi_2=f_1\circ\psi_1 :  \partial_+R_2 \lra \partial_+M
$$
preserves {\em oriented} contact structures (but the orientations of the characteristic foliations do not match). Since $f_s$ preserves the Legendrian curves $\psi_1(\gamma_\pm)$, the composition $\psi_2=f_1\circ\psi_1$ preserves the homotopy class of the intersection line field of $\DD^{(2)}$.
 
We will use the embedding $\psi_2$ in the proof of the following theorem. 
\begin{thm} \mlabel{t:M+N}
Let $M_1,M_2$ be manifolds with Engel structures $\DD_1,\DD_2$ such that both characteristic foliations admit closed transversals. Then $M_1\#M_2\#(S^2\x S^2)$ carries an Engel structure which coincides with $\DD_1$, respectively $\DD_2$, away from a neighborhood of the transversals where all connected sums are performed. The characteristic foliation of the new Engel structure again admits a closed transversal.
\end{thm}
\begin{proof}
Let us assume for the moment that the manifolds are oriented. Then the characteristic foliations are canonically oriented. For $i=1,2$ we cut $M_i$ along a closed transversal $N_i$ of the characteristic foliation. The resulting manifolds are still denoted by $M_1$ and $M_2$. The boundary of $M_i, i=1,2$ has two connected components $\partial_+M_i\simeq N_i\simeq\partial_-M_i$ and there is a natural identification
$$
\varphi : \partial_+M_1\cup \partial_+M_2 \lra \partial_-M_1\cup \partial_-M_2
$$
preserving oriented contact structures and the intersection line fields together with their orientations if $\DD_1,\DD_2$ are both oriented.

Let $U_1\subset\partial_+M_1$ and $U_2\subset\partial_+M_2$ be two balls. We will only need to modify intersection line fields on $U_1,U_2$ and $\varphi(U_1),\varphi(U_2)$. Therefore it is sufficient to orient $\DD_1,\DD_2$ only along $U_1,U_2$. From this we obtain coherent orientations of the intersection line fields on $U_1,\varphi(U_1)$ and $U_2,\varphi(U_2)$. 

In $U_1,U_2$ we choose Legendrian unknots $K_1,K_2$ with Thurston--Bennequin invariant $-2$ and rotation number $-1$. One can obtain $K_1,K_2$ by negative stabilization of the standard Legendrian unknot with Thurston--Bennequin invariant $-1$, cf.~\eqref{e:rot stab}.

Let $\psi_1$ be an embedding of $\partial_-R_1$ into $\partial_+M_1\cup\partial_+M_2$ which preserves oriented contact structures such that $\psi_1$ maps $\gamma_+$ to $K_1$ and $\gamma_-$ to $K_2$. Because the rotation numbers along $\gamma_\pm$ and $K_1, K_2$ equal $-1$, we can apply vertical modification to $\partial_+M_1\cup\partial_+M_2$ such that $\psi_1$ preserves oriented intersection line fields. 

Using the discussion above we obtain an attaching map 
$$
\psi_2 : \partial_+R_2 \lra \partial_-M_1\cup\partial_-M_2
$$
for $R_2$ with the model Engel structure $\DD^{(2)}$ which preserves oriented contact structures and the homotopy type of the intersection line field. Here we do not stick to the convention that $R_2$ is attached using an embedding of $\partial_-R_2$ but notice that the orientation of the characteristic foliations match. After a suitable vertical modification of the boundary, we obtain an Engel structure on 
$$
M=\left(\left(M_1\cup M_2\right)\cup_{\psi_1} R_1\right)\cup_{\psi_2} R_2~.
$$
Next we want to identify $\partial_+M$ with $\partial_-M$ such that oriented contact structures and the homotopy class of the intersection line fields are preserved. Away from $U_1\cup U_2$ we can use $\varphi$. This map has an obvious extension to a map $\partial_-M \lra \partial_+M$ which coincides with the natural identification of $\partial_+R_1\subset\partial_+M$ with $\partial_-R_2\subset\partial_-M$.

Unfortunately the contact structures on $\partial_+R_1$ and $\partial_-R_2$ are not preserved by this identification. However, the family of distributions defined by \eqref{e:alpha rho} away from $\partial_+R_1$ together with the restriction of \eqref{e:alpha s} to $\partial_+R_1\simeq \partial_-R_2$, is a family of contact structures. By \thmref{t:gray} there is an identification 
$$
\widetilde{\varphi} : \partial_-M \lra \partial_+M
$$
which preserves oriented contact structures and which coincides with $\varphi$ away from a neighborhood of the attaching region of $R_1$. 

Next we compare the intersection line fields on $\partial_+M$ and $\partial_-M$. Near $T^2_0=\{x=0\}\subset\partial_+R_1$ the Gray isotopy constructed in the proof of \thmref{t:gray} is induced by the Legendrian vector field 
\begin{equation} \label{e:Z}
Z=\frac{8y_2}{1+y_2^2}\left(y_2\rve{t}+\frac{1}{2}y_1\rve{x}\right)~.
\end{equation}
From \eqref{e:Z} one sees that $\widetilde{\varphi}$ maps the Legendrian curves $\{x=y_1=0\}\subset\partial_+R_1$ to themselves. The rotation number along these curves is $-1$ for both model Engel structures $\DD^{(1)}$ and $\DD^{(2)}$.

The rotation number of the intersection line field with respect to the Legendrian line field $\LL$ spanned by $y_2\partial_t+1/2\,y_1\partial_x$ is zero for both $\DD^{(1)}$ and $\DD^{(2)}$ (cf.~\propref{p:1 modelle}). Since $Z$ is tangent to $\LL$, the flow of this vector field preserves $\LL$. This implies that the identification $\widetilde{\varphi}$ preserves the homotopy class of the intersection line field along $T^2_0$.

Now let $\gamma$ be a curve in $\partial_+M$ which intersects $T^2_0\subset\partial_+R_1$ exactly once such that the endpoints of $\gamma$ lie outside of $U_1$ and $U_2$. We have to ensure that the intersection line fields are homotopic along $\gamma$ relative to the endpoints. For this we use the flexibility mentioned in \remref{r:suitable vert mod 1}: Choosing suitable values of $k$ for the vertical modification of the boundary when we attach $R_2$, we can ensure that $\widetilde{\varphi}$ preserves the homotopy class of the intersection line field along $\gamma$ relative to the endpoints. Then $\widetilde{\varphi}$ preserves the homotopy class of the intersection line field. 

After a vertical modification of the boundary we can identify $\partial_+M$ with $\partial_-M$ using $\widetilde{\varphi}$ such that we obtain an Engel structure on the resulting manifold $\widetilde{M}$. The new Engel structure coincides with $\DD_1$ and $\DD_2$ away from sufficiently big neighborhoods of $N_1$ and $N_2$. These transverse hypersurfaces are also contained in $\widetilde{M}$.

It remains to show that $\widetilde{M}$ is diffeomorphic to $M_1\# M_2\#(S^2\x S^2)$. In order to prove this, we decompose $R_1=h_1\cup h_2$ and $R_2=h_2'\cup R_3$, where $h_1,h_2,h_2',h_3$ are ordinary handles of index $1,2,2,3$, as in \eqref{e:decomp Rk}. The left part of \figref{b:unknotkirby} shows the attaching curve and the framing of $h_2$ after $h_1$ is attached. One end of $h_1$ lies in $U_1$ and the other end lies in $U_2$. 

Recall that $R_1$ and $R_2$ are attached to $M_1\cup M_2$ in a symmetric way. If we discard $h_2$ and $h_2'$ and identify the boundary components of $M_1\cup M_2$ with $h_1, h_3$ attached, then we obtain $M_1\# M_2$. 

Using an isotopy of the attaching curves of $h_2$ we can separate $h_2$, respectively $h_2'$, from $h_1$, respectively $h_3$, such that $h_2$ and $h_2'$ are attached along symmetric unknots in $\partial_+M_1\cup\partial_+M_2$, respectively $\partial_-M_1\cup\partial_-M_2$, which do not meet the attaching region of $h_1$, respectively $h_3$. 
\begin{figure}[htb] 
\begin{center}
\includegraphics[width=0.96\textwidth]{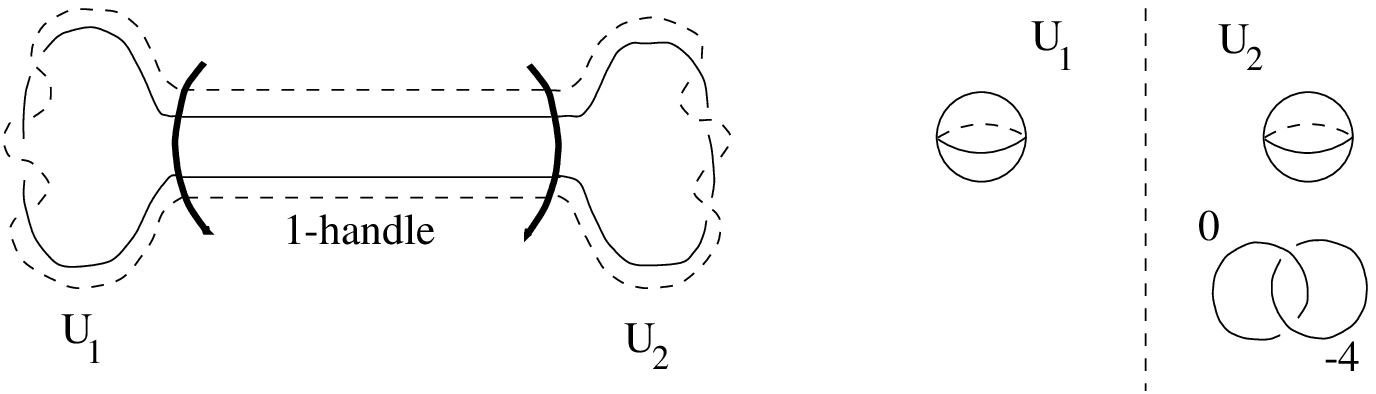}
\caption{\label{b:unknotkirby}}
\end{center}
\end{figure}

Since the Thurston--Bennequin invariant of $K_1$ and $K_2$ is $-2$,  the framing of the attaching curve of $h_2$ is $-4$ with respect to the obvious Seifert surface of the attaching curve. The $2$--handles $h_2, h_2'$ are attached by a doubling construction. The Kirby diagram therefore consists of spheres representing the attaching regions of $h_1$, the attaching curve of $h_2$ and a zero--framed meridian of the attaching curve of $h_2$, cf.~\cite{gompf}, p.~133. The ordinary handles of index $2$ therefore account for the summand $S^2\x S^2$. The right part of \figref{b:unknotkirby} shows the Kirby diagram of $M_1\cup M_2$ with $h_1, h_2, h_2'$ attached. 

This proves the claim under the assumption that $M_1$ and $M_2$ are oriented. If this is not the case we can nevertheless orient $\WW_i$ along the components of $\partial M_i$ after we have cut along the hypersurfaces. For the proof it suffices to orient $\WW_i$ along $U_i$ such that it points outwards and along $\varphi(U)$ such that it points inwards. All constructions carry over to this situation.
\end{proof}
In order to apply \thmref{t:M+N}, one has to find Engel structures whose characteristic foliation admits a closed transversal. This is true for the Engel structures from the proof of \thmref{T:EXIST}. We can also apply \thmref{t:M+N} to Engel structures obtained by prolongation after we perturb them slightly in the neighborhood of a leaf of the characteristic foliation (as in \propref{p:perturbed prolong}).
\begin{cor} \mlabel{c:PC1+PC2}
If $(N_1,\CC_1)$ and $(N_2,\CC_2)$ are contact manifolds, then $\PP\CC_1\#\PP\CC_2\#(S^2\x S^2)$ admits an Engel structure. 
\end{cor}
Using this corollary one obtains Engel structures on 
$$
M_n = n(S^3\x S^1)\#(n-1)(S^2\x S^2)~.
$$
One can show that it is impossible to construct an Engel structure on $M_n$ using the method of Geiges or prolongation, although $M_n$ is actually the total space of an orientable circle bundle over a $3$--manifold. However, the Euler class of all circle bundles with total space $M_n$ is odd (notice that since $M_n$ is orientable, the characteristic foliation of all Engel structures on $M_n$ are orientable). Since the tangent bundle of an orientable $3$--manifold is trivial, $M_n$ can not arise as circle bundle of a subbundle of rank $2$ of the $3$--manifold. 

Let us return to the proof of \thmref{t:M+N} and discuss the meaning of the assumption that both Engel structures have characteristic foliations which admit a closed transversal. We do not make explicit use of the fact that $N_1$ and $N_2$ are {\em closed} transversals in the proof. But implicitly, this assumption is used when we apply vertical modification of the boundary. If $N_1$ or $N_2$ have boundary, then we have to ensure that everything is well defined at boundary points of $N_i$ when we identify the boundary components of $M$. One can replace the existence of closed transversals in \thmref{t:M+N} by an assumption on the behavior of the Engel structure along leaves of the characteristic foliation.

Let $\WW(p)$ be a leaf of the characteristic foliation of an Engel structure $\DD$. We assume that $\WW(p)$ is not closed. Using nearby leaves of the characteristic foliation we obtain identifications of $\EE_p/\WW_p$ with $\EE_q/\WW_q$ for all $q$ on $\WW(p)$ by \lemref{l:transverse contact structure}. This allows us to define
\begin{align*}
\delta : \WW(p) & \lra \PP(\EE_p/\WW_p)\simeq S^1 \\
              q & \lmt [\DD_q]
\end{align*}
where we identify $\EE_q/\WW_q$ with $\EE_p/\WW_p$. This map is called {\em development map} of $\WW(p)$ and it is an immersion since $[\DD,\DD]=\EE$, cf.~\cite{adachi, montgomery}. For a given orientation of $\WW(p)$ let $\WW^\pm(p)$ be the segments of $\WW$ which lie on the two sides of $p$. Now we define the {\em twisting number}
$$
\tw^+(\WW(p))=\left|\delta^{-1}([\DD_p])\cap\WW^+(p)\right|~,
$$
i.e. $\tw^+(p)$ is the number of half twists of $\DD$ around $\WW$ in $\EE$ when one moves from $p$ in the sense of the orientation along the leaf. We define $\tw^-(\WW(p))$ similarly. Finally, if $\WW(p)$ is closed we define $\tw(\WW(p))$ as the degree of $\delta : \WW(p)=S^1 \lra S^1$. Now we can state a modified version of \thmref{t:M+N}.
\begin{thm} \mlabel{t:M+N zweite vers.}
Let $M_1,M_2$ carry Engel structures $\DD_1,\DD_2$ such that the characteristic foliation has non--closed leaves $\WW_1(p_1)$ through $p_1\in M_1$ and $\WW_2(p_2)$ through $p_2\in M_2$ with the property
\begin{equation} \label{e:twist bedingung M+N}
\tw^\pm(\WW_1(p_1))\ge C \textrm{ and } \tw^\pm(\WW_2(p_2))\ge C
\end{equation}
for some positive constant $C$ which is independent of the Engel structures.  

Then there is an Engel structure $\DD$ on $M_1\# M_2\#(S^2\x S^2)$ which coincides with $\DD_1,\DD_2$ outside of neighborhoods of $p_1\in\WW_1$ and $p_2\in\WW_2$. Moreover one can choose $\DD$ such that there is a leaf of the characteristic foliation which satisfies \eqref{e:twist bedingung M+N}. If $\WW_i(p_i)$ is closed one can replace \eqref{e:twist bedingung M+N} by $\tw(\WW_i(p_i))\ge 2(C+1)$ for $i=1,2$.  
\end{thm}
We just give a sketch of the proof and omit the details. Let $p_1,p_2$ satisfy the assumption of the theorem. For $i=1,2$ we cut $M_i$ along an embedded closed $3$--ball $N_i$ transverse to $\WW_i$, which is so small that the twisting numbers of all leaves of the characteristic foliation satisfy \eqref{e:twist bedingung M+N} after we have cut out $N_i$. Away from the boundary of $N_i$, the resulting space is a manifold with boundary. We can perform all constructions from the proof of \thmref{t:M+N} as long as we do not change anything near $\partial N_i$. 

Vertical modifications of the boundary are the only operations in the proof of \thmref{t:M+N} which affect $\partial N_i$. The problematic step is the vertical modification of the boundary used for the identification of the boundary components of $M$. The vertical modification of the boundary needed at this point is not allowed to change anything near $\partial N_i$. This means that the function $\widehat{f}$ used in \secref{ss:vert mod} for the description of a vertical modification of the boundary has to be zero near boundary points of $N_i$. The condition \eqref{e:twist bedingung M+N} allows us to choose $\widehat{f}\equiv 0$ on a collar of the boundary of $N_i$.  

We did not determine bounds for $C$. Since there are Engel structures (e.g. the standard Engel structure on $\R^4$) where the rotation number along all leaves of the characteristic foliation is zero, condition \eqref{e:twist bedingung M+N} is not fulfilled in general.

%------------------------ Bibliography

\end{document}